\documentclass[11pt]{article}
\usepackage{geometry}
\usepackage{natbib}
 \bibpunct[, ]{(}{)}{,}{a}{}{,}%

 \usepackage{endnotes}
\let\footnote=\endnote
\let\enotesize=\normalsize
\def\notesname{Endnotes}%
\def\makeenmark{$^{\theenmark}$}
\def\enoteformat{\rightskip0pt\leftskip0pt\parindent=1.75em
  \leavevmode\llap{\theenmark.\enskip}}

\usepackage{natbib}
 \bibpunct[, ]{(}{)}{,}{a}{}{,}%
 \usepackage{bm}
 \usepackage{bbm}
 \usepackage{amsmath}
 \usepackage{amssymb}
\usepackage{amsfonts}
\usepackage{graphicx}
\usepackage{subfigure}
\usepackage{float}
\usepackage{longtable,booktabs}
\usepackage{algorithm}
\usepackage{algorithmicx}
\usepackage{algpseudocode}
\usepackage{multirow}
\usepackage{epstopdf}
\usepackage[colorlinks,
            linkcolor=blue,
            anchorcolor=blue,
            citecolor=blue,
            ]{hyperref}
\usepackage{multirow}
\newenvironment{myproof}{{\indent  \it Proof:~}}{\hfill $\blacksquare$\par}
\usepackage[colorlinks,
            linkcolor=blue,
            anchorcolor=blue,
            citecolor=blue,
            ]{hyperref}

\usepackage[capitalise,noabbrev]{cleveref}
\usepackage{epstopdf}
\def\R{{\mathbb{R}}}

\def\e{{\mathbf e}}

\def\r{{\mathbf r}}

\def\v{{\mathbf v}}
\def\w{{\mathbf w}}
\def\x{{\mathbf x}}
\def\y{{\mathbf y}}
\def\z{{\mathbf z}}
\def\bz{{\mathbf 0}}
\def\cX{{\mathcal X}}
\def\Card{{\rm Card}}
\def\conv{{\rm conv}}
\def\argmin{{\mathop{\text{argmin}}}}
\def\blambda{{\boldsymbol\lambda}}
\usepackage{soul}

\newtheorem{fact}{Fact}
\newtheorem{procedure}{Procedure}
\usepackage{soul}
\newtheorem{theorem}{Theorem}
\newtheorem{lemma}{Lemma}
\newtheorem{definition}{Definition}
\newtheorem{remark}{Remark}
\newtheorem{assumption}{Assumption}

\begin{document}




\title{Decision Making under Cumulative Prospect Theory: An Alternating Direction Method of Multipliers}

\author{Xiangyu Cui
\thanks{School of Statistics and Management, Shanghai University of Finance and Economics, Shanghai, China, {cui.xiangyu@mail.shufe.edu.cn}}
\and
        Rujun Jiang
        \thanks{Corresponding author. School of Data Science, Fudan University, Shanghai, China, {rjjiang@fudan.edu.cn}}
        \and
        Yun Shi
       	\thanks{Academy of Statistics and Interdisciplinary Sciences, Faculty of Economics and Management, East China Normal University, Shanghai, China, {yshi@fem.ecnu.edu.cn}}
               \and
               Rufeng Xiao
        \thanks{School of Data Science, Fudan University, Shanghai, China, {rfxiao21@m.fudan.edu.cn}}
         \and
        Yifan Yan
        \thanks{School of Data Science, Fudan University, Shanghai, China,{yanyf21@m.fudan.edu.cn}}
}

\maketitle
\abstract{
This paper proposes a novel numerical method for solving the problem of decision making under cumulative prospect theory (CPT), where the goal is to maximize utility subject to practical constraints, assuming only finite realizations of the associated distribution are available. Existing methods for CPT optimization rely on particular assumptions that may not hold in practice. To overcome this limitation, we present the first numerical method with a theoretical guarantee for solving CPT optimization using an alternating direction method of multipliers (ADMM).  One of its subproblems involves optimization with the CPT utility subject to a chain constraint, which presents a significant challenge. To address this, we develop two methods for solving this subproblem. The first method uses dynamic programming, while the second method is a modified version of the pooling-adjacent-violators algorithm that incorporates the CPT utility function. Moreover, we prove the theoretical convergence of our proposed ADMM method and the two subproblem-solving methods. Finally, we conduct numerical experiments to validate our proposed approach and demonstrate how CPT's parameters influence investor behavior using real-world data.

}%



\maketitle

%


\section{Introduction}
Cumulative Prospect Theory (CPT, \citealp{KT:79}, \citealp{tversky1992advances}) is a popular model for descriptive decisions under risk and uncertainty in behavioral economics.  It describes four behaviour characteristics of the decision makers: (i) they assess the random outcomes as gains or losses with respect to a reference point, (ii) they are loss-averse, which means that they suffer more from losing than they enjoy from winning by the same amount, (iii) they are risk-averse for gains and risk-seeking for losses, and (iv) they assign more subjective weights to the outcomes with low objective probability. The CPT model  is widely applied in various fields, such as portfolio selection (\citealp{he2011portfolio}), operations management (\citealp{nagarajan2014prospect}), supply chain management (\citealp{liu2013order}), transportation management (\citealp{li2015multi}), etc.

The CPT model incorporates an S-shaped utility function defined as follows:
	\begin{equation}\label{eq:sshapeu}
	U(z)=\left\{\begin{array}{ll}   -\mu (B-z)^{\alpha} ,& z \leq B,  \\ (z-B)^{\alpha} , & z > B,\end{array}\right.
	\end{equation}
where $B$ serves as the reference point distinguishing gains ($z > B$) from losses ($z \leq B$). The parameter $\mu > 1$ signifies greater distress from losses compared to joy from equivalent gains, known as the loss aversion parameter. For gains ($z > B$), the utility function is concave, indicating risk aversion, while for losses ($z \leq B$), it is convex, implying risk-seeking behavior. The parameter $\alpha$ reflects the degree of relative risk aversion (or seeking) concerning gains (or losses). Additionally, the CPT model integrates inverse S-shaped probability weighting functions for gains and losses:
	\begin{equation}\label{eq:pd}	
\begin{array}{rll}
	 \omega_-(F_Z(z)) &=~ \frac{F_Z(z)^{\delta}}{((1-F_Z(z))^{\delta}+ F_Z(z)^{\delta})^{1/\delta}}, &\quad z\leq B,\\
	\omega_+(1-F_Z(z)) &= ~\frac{(1-F_Z(z))^{\gamma}}{((1-F_Z(z))^{\gamma}+ F_Z(z)^{\gamma})^{1/\gamma}},&\quad z>B,
	\end{array}
	\end{equation}
where $F_Z(z)$ represents the cumulative distribution function (CDF) of the random outcome $Z$. Here, $\delta$ and $\gamma$ are the probability distortion parameters for losses and gains, respectively. In the loss domain, $\omega_-(F_Z(z))$ distorts $F_Z(z)$, with $\delta < 1$ indicating that the decision maker assigns higher subjective probabilities to small outcomes than their objective probabilities. In the gain domain, $\omega_+(1-F_Z(z))$ distorts the survival function $1-F_Z(z)$, where $\gamma < 1$ implies assigning higher subjective probabilities to large outcomes than their objective probabilities. The decision problem faced under the CPT model is represented as:
\begin{equation}\label{eq:continue-1}
\begin{aligned}
\max_{\x\in\R^d} ~\int_{B}^{+\infty} (z-B)^{\alpha}  d [-\omega_+ (1-F_Z(z))] -\mu \int_{-\infty}^{B} (B-z)^{\alpha}  d \omega_- (F_Z(z))\quad {\rm s.t.}~  z= \r^T \x,\quad \x \in \mathcal X,
\end{aligned}
\end{equation}
where the returns of $d$ assets comprise a vector $\r$, which is a $d$-dimensional random vector, the portfolio weight $\x \in \mathbb{R}^d$  is the decision variable, $^T$ denotes the transpose operator, and  $\mathcal X$ is the feasible set for decision variable $\x$. The portfolio return $z$ is a scalar and is a random variable whose randomness comes from $\r$. It is important to note that both the utility function and the probability weighting functions exhibit nonconvexity and nonsmoothness at the reference point, posing challenges in solving decision-making problems within the CPT framework.

By using the notations of the general utility and probability weighting functions,
the decision problem faced under the CPT model can be rewritten in a more compact form as follows,
\begin{equation}\label{eq:continue-2}
\begin{aligned}
\max_{\x\in\R^d}  ~ \int_{B}^{+\infty} U(z) d [-\omega_+ (1-F_Z(z))] + \int_{-\infty}^{B}  U(z) d \omega_- (F_Z(z))\quad {\rm s.t.}~  z= \r^T \x,\quad \x \in \mathcal X.
\end{aligned}
\end{equation}
Let $q_\nu$ be the $\nu$th quantile of $z$, where $\nu \in (0,1)$.  The objective function in \eqref{eq:continue-2} has an equivalent form.
 \begin{equation}\label{eq:dist5}
    \begin{aligned}
        & \int_{B}^{+\infty} U(z) d [-\omega_+ (1-F_Z(z))] + \int_{-\infty}^{B}  U(z) d \omega_- (F_Z(z))\\
        = & \lim_{N\to +\infty} \sum_{i=1}^N U(q_\frac{i}{N}) \left\{ \left[ \omega_- \left(\frac{i}{N}\right) - \omega_- \left( \frac{i-1}{N}\right)\right] \mathbbm1_{\{ q_{\frac{i}{N}}\leq B\}} \right.\\
        & \quad \quad \left.+ \left[ \omega_+ \left(1 - \frac{i-1}{N}\right) - \omega_+ \left(1 - \frac{i}{N}\right)\right] \mathbbm1_{\{ q_{\frac{i}{N}}> B\}} \right\},
\end{aligned}
\end{equation}
where $\mathbbm1_{\{\cdot\}}$ is the characteristic function.
In real-world applications, it is difficult to know the true distribution of the random variable $\r$. Instead, we can often access finite realizations of the distribution via either historical data or sample points generated by simulation. The estimated distribution of $z$ is determined by the realizations of  $\r$.
Let $R\in \R^{N\times d}$ be the matrix formed by $N$ realizations of $\hat{\r}_i,i =1,2, \cdots,N$.
The estimated $\frac{i}{N}$th quantile of $z$ is $(R\x)_{[i]}$, where we use subscript $[i]$ to denote the $i$th smallest element of a vector. Define a piecewise function
\begin{equation}\label{eq:defc}
c_i(q_{\frac{i}{N}}):= \left\{ \left[ \omega_- \left(\frac{i}{N}\right) - \omega_- \left( \frac{i-1}{N}\right)\right] \mathbbm1_{\{ q_{\frac{i}{N}}\leq B\}} + \left[ \omega_+ \left(1 - \frac{i-1}{N}\right) - \omega_+ \left(1 - \frac{i}{N}\right)\right] \mathbbm1_{\{ q_{\frac{i}{N}}> B\}} \right\}.
\end{equation}
This equation indicates that $c_i$ takes on two distinct values on either side of  $B$.
For simplicity, we will use the convention $c_i=c_i(q_{\frac{i}{N}})$.
Then we derive the discrete form of CPT optimization
\begin{equation} \label{pb:1}
\begin{aligned} \min_{\x\in\R^d} &~ -\sum_{i = 1}^N c_i U((R\x)_{[i]})
  \quad &{\rm s.t.}~\x \in \mathcal X. \end{aligned}
\end{equation}
 The difficulty of this problem stems from the possible nonconvexity and nonsmoothness of the objective function. Indeed,  the S-shape utility function $U$ exhibits nonsmoothness at the reference point and is nonconvex, the rank of the historical return introduces further nonsmoothness, and the piecewise nature of the CPT weight $c_i$  adds nonconvexity and nonsmoothness as well.
To the best of our knowledge, the problem does not fall into the scope of any existed non-convex or non-smooth solvers. It cannot be tackled by subgradient-type methods due to its nonconvexity \citep{boyd2003subgradient}, and it lacks the special composite structures needed for the application of proximal gradient-type methods \citep{parikh2014proximal}.
This realization motivates us to develop a practical method tailored to solving problem (7) by leveraging its particular structure.

\subsection{Related Work}
A significant portion of the literature addresses these challenges by introducing specific assumptions that simplify problem structures but may not accurately represent real-world scenarios. For example, \cite{BH:08,liu2013order} propose the normal distribution assumption, which is used to model the portfolio's outcome $Z$ according to a normal distribution, facilitating easier comparisons between outcomes. This approach simplifies the comparison of portfolio outcomes by examining whether one outcome's CDF crosses another outcome's CDF from below, as demonstrated in Proposition A2 of \cite{BH:08}. Similarly, \cite{shi2015discrete} introduces the elliptical distribution assumption, while \cite{Barberis/Xiong:06} introduces the binomial tree assumption. Other works, such as \cite{bernard2010static,he2011portfolio,nagarajan2014prospect}, propose the one-dimensional random source assumption, which transforms the decision-making problem into a simpler one-dimensional optimization problem, even if the problem is not convex, making it manageable. \cite{jin2008behavioral,van2020dynamic} propose the complete market assumption, enabling the use of the martingale method and converting portfolio decision-making into portfolio's return decision-making. Additionally, \cite{De:2012,shi2015dynamic} propose the piece-wise linear utility function assumption. However, it is important to note that this assumption overlooks nuances such as risk aversion in gains and risk-seeking behavior in losses, resulting in a convex utility function that ignores the complexities inherent in real-world decision-making processes.

On the other hand, there are only limited numerical algorithms for solving the CPT models with a small number of assets, such as heuristics
based algorithms \cite{barro2020cumulative}, or
grid search methods \cite{hens2014cumulative}.
Very recently, \cite{luxenberg2024portfolio} propose three numerical optimization methods for solving CPT optimization. However, their methods either only apply to an approximation of the original CPT optimization problem by enforcing monotonicity on weights, or are heuristics without theoretical guarantees. To the best of our knowledge, there is no existing method that can solve the CPT optimization problem with a general setting.

\subsection{Main Contribution}

We aim to solve the CPT optimization problem with a general setting by the alternating direction method of multipliers (ADMM). The ADMM is a popular optimization method for problems that have a certain separable structure in  objective functions (\citealp{glowinski1977numerical}).
Its main advantage is that using the separable structure, the subproblems in ADMM are easy to solve.
The ADMM finds numerous applications in real world such as statistical learning \citep{scheinberg2010sparse}, portfolio optimization \citep{cui2018portfolio} and signal processing \citep{combettes2007douglas}.
For more details on the ADMM algorithm, we refer to the survey paper \cite{boyd2011distributed}.


The main contribution of this paper is to propose an ADMM algorithm.
 To address the difficulty in the utility function, we introduce an auxiliary variable $\y\in\R^N$ to represent return realizations of the portfolio. We first reformulate \eqref{pb:1} by introducing an auxiliary variable  $\y=R\x$,
\begin{equation} \label{pb:umy}
\begin{aligned}
\min_{\x,\y} &~-\sum_{i = 1}^N c_i U(y_{[i]})+I_{\mathcal X}(\x)
&\quad {\rm s.t.} ~\y = R\x,
\end{aligned}
\end{equation}
where $I_{\mathcal X}(\x)=0$ if $\x\in\mathcal X$ and $\infty$ otherwise is the indicator function of set $\mathcal X$.
Utilizing the separability of $\x$ and $\y$ in the objective function of the new reformulation, we then develop an ADMM as a solution method and  demonstrate its convergence property.
 The primary advantage of ADMM is that it decomposes the original problem into subproblems that are easier to solve. In deed, the $\x$-subproblem of the ADMM (to be defined in \cref{sec:3})
$$\min_{\x\in \mathcal X}~ \frac{\sigma}{2} \left \| \y^k - R\x + \frac{\blambda^k}{\sigma} \right \|_2^2$$
is a convex program with a convex quadratic objective function provided  $\cX$ is a convex set, which is easy to solve.
The $\y$-subproblem (to be defined in \cref{sec:3}),  involving a nonconvex finite sum objective with rank-dependent weights, can be reformulated as a chain-constrained separable problem with a favorable structure
\begin{equation*}
        \min_{\y} ~ \sum_{i=1}^N -c_iU(y_{l_i}) + \frac{\sigma}{2}(y_{l_i} - w_{l_i})^2  \quad
        \text{s.t.} \quad y_{l_1} \leq y_{l_2} \leq \dots\le y_{l_N}.
\end{equation*}
See paragraph before \eqref{eq:yisot} in \cref{sec:4} for the derivation.
Albeit challenging to solve, the advantageous structure of the $\y$-subproblem inspires us to devise two customized algorithms for its solution.
The first one is a dynamic programming (DP) method. We demonstrate that the DP method can find the global optimal solution. Despite the solidness in theory, the DP method may be slow in practice if there are a large number of historical scenarios of the randomness.
To improve practical efficiency, we further propose a variant of the celebrated pooling-adjacent-violators (PAV) algorithm to solve this subproblem, which can be seen as an extension of the PAV algorithm for convex programming \citep{ayer1955empirical,brunk1972statistical,best1990active,ahuja2001fast}.
We further prove that the PAV algorithm converges to a stationary point of a reformulation of the subproblem.
Our numerical study shows that the proposed ADMM outperforms the general solver \texttt{fmincon} in MATLAB { and methods proposed in \cite{luxenberg2024portfolio}}, and the ADMM with PAV as its subproblem performs better than that equipped with the DP.
We also conduct an empirical portfolio selection study based on 48 industry indices of the US market  and investigate the impacts of CPT's parameters on the optimal portfolio.

The remainder of this paper is organized as follows.
In \cref{sec:2}, we specify the expression of the decision making problem under CPT. Then in \cref{sec:3}, we propose the ADMM framework and derive its convergence analysis. In \cref{sec:4}, we propose two methods, the DP and the PAV methods, to solve the $\y$-subproblem in the ADMM. In \cref{sec:5}, we conduct numerical experiments to show the effectiveness of the proposed method.
We conclude our paper in \cref{sec:7}.

\section{The decision making problem under CPT}\label{sec:2}
In this paper, we investigate the decision making problem under CPT \eqref{pb:1}.

We point out that the CPT optimization problem  \eqref{pb:1} can represent a portfolio optimization problem whose objective is to maximize the CPT utility, where $R$ represents the matrix formed by $N$ historical returns of $d$ assets,  $R\x$ represents the historical return of the portfolio $\x$, and $c_i$ is the weight associated with the $i$th order statistics of the historical return of the portfolio.

 This formulation also includes many financial optimization applications as its special cases. When we choose a strictly concave utility and a single inverse S-shaped probability weighting function, the portfolio optimization model under CPT becomes the rank dependent utility maxization problem (RDU, \citealp{quiggin2012generalized}). Furthermore,
when we choose the identity map $U(y)=y$ and different weights $c_i$, problem \eqref{pb:1} includes many risk management problems. When $c_i=1$ for $i=\lceil \alpha N \rceil$ and $c_i=0$ otherwise,
the objective function reduces to Value at Risk at confidence level $\alpha\in(0.5,1)$ (VaR, \citealp{duffie1997overview},\citealp{cui2018portfolio}), where $\lceil \alpha N \rceil$ is the smallest integer that is larger than or equal to $ \alpha N $.
When $c_i=\frac{1}{N-\lceil \alpha N \rceil+1}$ for $i\ge\lceil \alpha N \rceil$ and $c_i=0$ otherwise,
the objective function reduces to Conditional Value at Risk at confidence level $\alpha$ (CVaR, \citealp{rockafellar2000optimization}).
When $\sum_{i=1}^{N}c_i=1$ and $0\le c_1\le \dots\le c_N\le 1$, then the objective function is the spectral risk measure (SRM, \citealp{acerbi2002spectral}). When $\sum_{i=1}^{N}c_i=1$ and $0\le c_i$, then the objective function is the distortion risk measure (DRM, \citealp{dhaene2006risk}). As the { rank-dependent} weights $c_i$ in VaR and distortion risk measure are not decreasing functions of $i$, the corresponding risk management problems are not convex and hard to solve.

To make our research applicable to a more general setting, we consider the decision making problem under CPT \eqref{pb:1}, which does not specify the concrete form of the utility function, and only requires the utility function { to satisfy} the following assumption.
\begin{assumption}\label{asmp:1}
Assume the function $U:\R\to\R$ in problem \eqref{pb:1} satisfies the following properties:
1. $U$ is a strictly increasing, continuous function in $\R$.
2. $U$ is a third order continuously differentiable function in either $(-\infty,B)$ or $(B,\infty)$.
3. When $z \in (-\infty, B)$, we have $U''(z)>0$. When $z \in (B, \infty)$,  we have $U''(z)<0$.
4.  When $z \in (-\infty, B)$,  we have $U{'''}(z) > 0$.
5.  If $U{'}(B^{-}) = \infty$, then $U{'}(B^{+}) = \infty$, and vice versa,  where $U{'}(\theta^+) = \lim_{y\rightarrow\theta^+}U{'}(y)$, and $U{'}(\theta^-) = \lim_{y\rightarrow\theta^-}U{'}(y)$.
\end{assumption}
 The first property represents the fundamental necessity for utility functions to accurately reflect preferences (\citealp{von1947theory}). The second property constitutes a straightforward requirement for the utility function to exhibit smoothness (\citealp{eeckhoudt2006putting}). The third property specifies the essential condition for the CPT model, indicating risk aversion in the gain domain and risk-seeking behavior in the loss domain (\citealp{tversky1992advances}). The fourth property explores the nuanced nature of convexity tendencies in the loss domain, particularly emphasizing their strengthening near the reference point. This property elucidates that the marginal utility in the loss domain is a convex function, which represents the decision maker is prudent (\citealp{eeckhoudt2006putting}). The fifth property underscores the need for utility functions in both gain and loss domains to demonstrate similar first-order behaviors, ensuring consistency in decision-making processes. The first, second, fourth properties serve as relatively mild prerequisites for utility functions. In contrast, the third property is crucial for defining feature (iii) within the CPT model. The last property is designed for power utility. When the first-order derivatives around the reference point are finite, the fifth property is not needed.

The utility function proposed by \cite{tversky1992advances} satisfies Assumption \ref{asmp:1}, and the outlined properties extend to a generalized version of this utility function. For example, consider the following utility function:
\begin{equation*}
V(z)=\left\{\begin{array}{ll} \mu (e^{\alpha_- (z-B)}-1),& z \leq B,  \\  1-e^{-\alpha_+ (z-B)}, & z > B,\end{array}\right.
\end{equation*}
where $\alpha_+>0$, $\alpha_->0$, which also meets Assumption \ref{asmp:1}. When $\mu=1$ and $B=0$, it reduces into the utility function in \cite{luxenberg2024portfolio}.  Following \cite{pratt1964risk}, we can compute the absolute risk aversion coefficient in the gain domain and absolute risk-seeking coefficient in the loss domain as follows:
\begin{equation*}
-\frac{V''(z)}{V'(z)}=\left\{\begin{array}{ll} -\alpha_-,& z \leq B,  \\  \alpha_+, & z > B,\end{array}\right.
\end{equation*}
indicating constant absolute risk aversion and risk-seeking coefficients in this utility function. Notably, this utility function differs from the power utility-based function proposed by \cite{tversky1992advances}. For the utility function in \eqref{eq:sshapeu}, we can compute the absolute risk aversion coefficient in the gain domain and absolute risk-seeking coefficient in the loss domain as follows:
\begin{equation*}
-\frac{U''(z)}{U'(z)}=\left\{\begin{array}{ll} (1-\alpha)(B-z)^{-1},& z \leq B,  \\ (1-\alpha)(z-B)^{-1}, & z > B,\end{array}\right.
\end{equation*}
indicating decreasing absolute risk aversion and increasing absolute risk-seeking coefficients.

Assumption \ref{asmp:1} does not impose specific constraints on probability weighting functions, allowing for flexibility in incorporating various functions proposed in the literature. For instance, the model can accommodate probability weighting functions outlined in \cite{tversky1995weighing} as well as those in \cite{prelec1998probability}. The probability weighting functions in \cite{tversky1995weighing} are expressed as:
\begin{align*}
&\omega_+(1-F_Z(z)) = \frac{\gamma_+ (1-F_Z(z))^{\delta_+} }{\gamma_+ (1-F_Z(z))^{\delta_+}  + (F_Z(z))^{\delta_+}},\quad \omega_- (F_Z(z)) = \frac{\gamma_- (F_Z(z))^{\delta_-} }{\gamma_- (F_Z(z))^{\delta_-}  + (1-F_Z(z))^{\delta_-}},
\end{align*}
where $0<\delta_+, \delta_-<1$ and $\gamma_+, \gamma_->0$ are parameters that influence the weighting behavior.  Similarly, \cite{prelec1998probability} proposes probability weighting functions given by:
\begin{align*}
& \omega_+(1-F_Z(z)) = e^{-\gamma_+ (-\ln(1-F_Z(z)))^{\delta}},\quad  \omega_-(F_Z(z)) = e^{-\gamma_- (-\ln(F_Z(z)))^{\delta}},
\end{align*}
where $0<\delta<1$ and $\gamma_+, \gamma_->0$ are parameters affecting the weighting patterns. Moreover, the model can accommodate scenarios where decision makers either overestimate or underestimate the probabilities of extreme outcomes. This phenomenon, documented in studies like \cite{polkovnichenko2013probability} and \cite{shi2023beta}, can be captured by setting the probability distortion parameters to values less than 1 (overestimation) or greater than 1 (underestimation), respectively.

\section{The ADMM Algorithm}\label{sec:3}
In this section, we propose an alternating direction method of multipliers (ADMM) to solve problem \eqref{pb:1}. Among all the utilities used in the literature, the CPT utility poses the most difficulty for optimization.
Therefore, we focus on the CPT optimization in this section and give a remark on other risk preferences that our method can handle at the end of this section.
Note that the main challenge of optimizing \eqref{pb:1} stems from the utility function, which is nonsmooth and nonconvex.
At first glance, the nonsmoothness arises from the order statistics of $R\x$, which are finite realizations of the return of the portfolio, and the nonconvexity arises from the S-shaped utility. However, we point out that a more subtle issue is that the weight $c_i$ in \eqref{eq:defc}, which is determined by the probability distortion, the reference point $B$ and rank of realizations of the return of the portfolio, introduces additional nonconvexity and nonsmoothness.
The non-convexity and the non-smoothness made it difficult to directly employ the existing numerical solvers.


This motivates us to consider the reformulation \eqref{pb:umy} that separates the constraint of the portfolio strategy and the utility function with the auxiliary variable $\y=R\x$.
The augmented Lagrangian function of   problem \eqref{pb:umy} is
\begin{equation*}\label{eq:alm}
        L_{\sigma}(\x,\y;\blambda) = - \sum_{i = 1}^{N} c_i U(y_{[i]}) + I_{\mathcal X}(\x) + \langle\blambda,\y  - R\x\rangle +\frac{\sigma}{2} \left \|  \y - R\x   \right \| _2^2,
\end{equation*}
where $\blambda$ is the Lagrange multiplier and $\sigma > 0$ is the quadratic penalty parameter.
Note that it is difficult to simultaneously optimize $\x$ and $\y$ in the augmented Lagrangian function.
This motivates us to adopt the ADMM framework that is widely used to handle a separable objective function.
The update of  ADMM is as follows
\begin{subequations}
    \begin{align}
     \x^{k+1} &= \argmin_{\x \in \mathcal X }~L_{\sigma}(\x,\y^k;\blambda^k),\label{eq:subx} \\
     \y^{k+1} &= \argmin_{\y} ~L_{\sigma}(\x^{k+1},\y;\blambda ^k),\label{eq:suby} \\
     \blambda^{k+1} &= \blambda^{k} + \sigma( \y^{k+1} - R\x^{k+1} ).\label{eq:sublam}
    \end{align}
\end{subequations}
We summarise our ADMM framework in \cref{alg:admm}.

The $\x$-subproblem is
$$\min_{\x} ~I_{\mathcal X}(\x) + \langle\blambda^k,\y^k  - R\x\rangle +\frac{\sigma}{2} \left \|  \y^k - R\x   \right \| _2^2,$$
which is independent of the utility function.

When optimizing the $\x$-subproblem, $\y^k$ and $\blambda^k$ are seen as constants, and thus the $\x$-subproblem is equivalent to
$$\min_{\x\in \mathcal X}~ \frac{\sigma}{2} \left \| \y^k - R\x + \frac{\blambda^k}{\sigma} \right \|_2^2.$$

The $\x$-subproblem \eqref{eq:subx} is  a convex program  when $\mathcal X$ is a convex set, which can be easily solved by off-the-shelf solvers.

\begin{center}
\begin{minipage}{\linewidth}
\begin{algorithm}[H]\small
        \caption{ADMM for CPT optimization problem \eqref{pb:1}}
        \label{alg:admm}
        \begin{algorithmic}[1]
        \Require  $R,~\sigma,~\epsilon_1,~\epsilon_2,~\y^0,~\blambda^0 $, and $k = 0$
            \Repeat

         \State Solve \eqref{eq:subx} to obtain $\x^{k+1}$;
            \State Solve \eqref{eq:suby} to obtain $\y^{k+1}$ either exactly or approximately;
            \State Update $\boldsymbol{\lambda^{k+1}} = \boldsymbol{\lambda^{k}} + \sigma( \y^{k+1} - R\x^{k+1} )$;
            $k = k + 1$
            \Until{$\| \y^{k} - R\x^{k} \| \le \epsilon_1$ and $\| \y^{k} - \y^{k-1}\| \le \epsilon_2$}
        \end{algorithmic}
\end{algorithm}
\end{minipage}
\end{center}

The main difficulty lies in solving the $\y$-subproblem \eqref{eq:suby}, which is a nonconvex and nonsmooth optimization problem.  Fortunately, by exploiting the properties of the utility function $U$ under \cref{asmp:1} and the piecewise structure of its weights $c_i$, we propose two algorithms to solve it, whose details are given in \cref{sec:4}.
At the $k$th iteration, the $\y$-subproblem is to minimize the following function
\begin{equation*}
        \Phi(\y) = - \sum_{i = 1}^{N} c_i U(y_{[i]}) + \frac{\sigma}{2} \left \|\y - R\x^{k+1} + \frac{\blambda^k}{\sigma} \right\|^2.
\end{equation*}

Before presenting our convergence analysis of the proposed ADMM algorithm, let us first introduce some results from nonsmooth analysis. Let $\Omega(\y)$ denote the first component of $\Phi(\y)$, i.e.,
$\Omega(\y) = - \sum_{i = 1}^{N} c_i U(y_{[i]}).$
We show that both $\Phi(\y)$ and $\Omega(\y)$ are locally Lipschitz functions under mild conditions so that the stationary point can be characterized by the Clarke generalized gradient \citep{clarke1990optimization}. 
Recall that $U’(B) < \infty$ represent that both $U’(B^-) < \infty$ and $U’(B^+) < \infty$.

\begin{lemma}\label{lem:locallip}
Under \cref{asmp:1}, the functions  $\Phi(\y)$ and $\Omega(\y)$ are both locally Lipschitz continuous in some neighbourhood of any given  $\bar\y$ satisfying $U’(\bar y_i) < \infty, \forall i=1,\ldots,N$.
\end{lemma}
The proofs in this paper are all deferred to Appendix \ref{ap:proof}.
We remark that $\Omega(y)$ is non-Lipschitz at $y_i=B$ for the S-shape power utility \eqref{eq:sshapeu}.
As we have shown that $\Omega(\y)$ is locally Lipschitz for every $\y$ as long as $U'(y_i)< \infty ~\forall i$, the Clarke generalized gradient of $\Omega(\y)$, denoted by $\partial \Omega(\y)$, exists if $U'(y_i)< \infty~\forall i$ \citep{clarke1990optimization}.
Let $\Xi(\x) = \Omega(R\x)$. By the chain rule, we have $$\partial \Xi(\x) \subset R^T \partial \Omega(R\x).$$
Denote by $N_{\mathcal X}(\bar\x)$  the normal cone of $\mathcal X$ at  $\bar\x\in\cX$, i.e., $N_{\mathcal X}(\bar\x) = \{ \z \in \mathbb{R}^d \mid \langle\z, \x - \bar\x\rangle \leq 0,~\forall \x \in \mathcal X\}.$
It is well known that $\partial I_\cX(\x)=N_\cX(\x)$; see, e.g., \cite[Example 3.5]{beck2017first}.
We then have
\[
\partial\left(\Xi(\x)+I_{\mathcal X}(\x)\right)\subset \partial \Xi(\x) + \partial I_{\mathcal X}(\x)= \partial \Xi(\x) + N_{\mathcal X}(\x),
\] where the first inclusion is due to Corollary 1 of \cite[Theorem 2.9.8]{clarke1990optimization}.
Suppose $\x^*$ is a local minimizer of \eqref{eq:suby}. Then from  \cite[Proposition 2.3.2]{clarke1990optimization} we have
$
\bz\in \partial \left(\Xi(\x^*)+I_{\mathcal X}(\x^*)\right).
$
The above facts imply the following necessary optimality condition
\begin{equation}\label{eq:necon}
        \bz\in R^T \partial \Omega(R \x^*) + N_{\mathcal X}(\x^*).
\end{equation}

Let $\y$ be such that $\Omega$ is locally Lipschitz near $\y$. Due to \cite[Theorem 2.5.1]{clarke1990optimization}, we have
\begin{equation*}
\partial\Omega(\y)=\conv\left(\{\lim_{i\to\infty}\nabla \Omega(\y^i):~\y^i\to\y,~\Omega\text{ is differentiable at }\y^i\}\right),
\end{equation*}
where $\conv(A)$ denotes the convex hull of set $A$. 
In fact, using the structure of $\Omega(\y)$, we can characterize the expression of  $\partial\Omega(\y)$ explicitly by considering the coordinates that are equivalent. Let $\{i_1,i_2,\ldots,i_{s^k}\}$ be a permutation of $\{1,\ldots,N\}$  and $\y$ be such that
\begin{equation}
\label{eq: y array}
y_{i_{s^1+1}}=\cdots=y_{i_{s^2}}<\cdots <y_{i_{s^{k-1}}+1}=\cdots=y_{i_{s^k}},
\end{equation}
where we use the  convention $s^1=0$ and $s^k=N$.
Note that except the points that are equivalent to $B$, the only non-differentiable points are those where at least two coordinates have the same value.

We then have
\begin{equation}
\label{eq:clarke}
\begin{array}{ll}
\partial\Omega(\y)=\conv\left(\left\{\v:\right.\right.~v_{i_j}\in -c_{i_j}\partial U(y_{i_j}),~
j=s^t+1,\ldots,s^{t+1},~t=1,\ldots,k,\\
\quad \{i_1,i_2,\ldots,i_{s^k}\} \text{ is a permutation of }\{1,\ldots,N\} \left.\left.\text{ that satisfies \eqref{eq: y array}.}\right\}\right).
\end{array}
\end{equation}
We remark that $\partial U(y_{i_j})$ is always a singleton when $y_{i_j}\neq B$, and $\partial U(B)=\conv\{U(B^-),U(B^+)\}$ is an interval.
Let us illustrate this using a simple example. If $\y=(y_1,y_2,y_3)$ with $y_2=y_3<y_1$, then we must have either $y_{[1]}=y_2,~y_{[2]}=y_3,~y_{[3]}=y_1$ or $y_{[1]}=y_3,~y_{[2]}=y_2,~y_{[3]}=y_1$,
and thus
\[
\partial\Omega(\y)=\conv\left(\{(-c_{3}U'(y_1),-c_1U'(y_2),-c_2U'(y_3))^T,
(-c_{3}U'(y_1),-c_2U'(y_2),-c_1U'(y_3))^T\}\right).
\]
We then directly have from \eqref{eq:clarke}  the following lemma  that is useful in our convergence analysis.
\begin{lemma}[Outer semicontinuity]\label{lem:hemic}
Let $\y$ be such that $\Omega$ is locally Lipschitz near $\y$.
Then for any sequences $\{\v^k\}$ and $\{\y^k\}$ such that $\v^k\to\v^*$, $\y^k\to\y^*$ and $\v^k\in \partial\Omega(\y^k)$, we have
$\v^*\in \partial\Omega(\y^*)$.
\end{lemma}

As the main theoretical result in this section, we will show that the proposed ADMM algorithm converges to a point satisfying the necessary optimality condition \eqref{eq:necon} under mild conditions.
To this end, let us first make an assumption on the sequence of Lagrange multiplier.
\begin{assumption}\label{asmp:2}
The sequence $\{ \blambda^k \}$ is bounded and satisfies the following condition
        \begin{equation*}\label{lambdaequation}
                \sum_{k = 1}^{\infty} \| \blambda^{k+1} - \blambda^k \|^2 < \infty.
        \end{equation*}
\end{assumption}
We remark that the convergence of nonconvex ADMM with two nonsmooth blocks { is} challenging without imposing assumptions like Assumption \ref{asmp:2}; see, e.g., \cite{wang2019global,lin2022alternating}. Additionally, it is worth noting that \Cref{asmp:2} is extensively utilized in the ADMM literature, as seen in works like \cite{xu2012alternating} and \cite{bai2021augmented}. A stronger version of this assumption is used in \cite{shen2014augmented}.

We also impose the following assumption on the solutions of the $\x$- and $\y$-subproblem.
\begin{assumption}\label{asmp:3}
The $\x$-subproblem is globally solved. The $\y$-subproblem is solved such that
\begin{equation}\label{eq:ydescent}
        L_{\sigma}(\x^{k+1},\y^k;\blambda^k) - L_{\sigma}(\x^{k+1},\y^{k+1};\blambda^k) \geq 0,\quad \text{and}
\end{equation}
\begin{equation} \label{ykfirst}
\bz \in \partial\Omega(\y^{k+1}) + \sigma \left(\y^{k+1} - R\x^{k+1} + \frac{\blambda^k}{\sigma}\right).
\end{equation}
\end{assumption}
We remark that the above assumption is quite mild.  The assumption that the $\x$-subproblem is globally solved is easy to satisfy if $\mathcal X$ is a convex set, because such a convex program can be solved globally by many methods such as interior point methods \citep{nocedal1999numerical}. Meanwhile, the $\y$-subproblem is solved to a stationary point $\y^{k+1}$ that has an objective value not exceeding that of the initial point  $\y^{k}$. We will propose a subproblem solver that can return a solution satisfying both \eqref{eq:ydescent} and \eqref{ykfirst} in \cref{sec:4}.

Now we are ready to present the main convergence result for the ADMM algorithm.
\begin{theorem}\label{thm:admm}
Suppose that $R^TR$ is positive definite, $\cX$ is a convex and compact set, $\epsilon_1=\epsilon_2=0$, and Assumptions \ref{asmp:2} and \ref{asmp:3} hold. Let $\{(\x^k,\y^k)\}$ be a sequence generated by the ADMM algorithm.
Then any accumulation point of the sequence $\{(\x^k,\y^k)\}$, denoted by $(\x^*,\y^*)$,  satisfies the following condition
\begin{equation*}
        \bz\in R^T \partial \Omega(\y^*) + N_{\mathcal X}(\x^*),  \quad       \y^* = R \x^*.
\end{equation*}
Hence $\x^*$ is a stationary point of \eqref{pb:1} in the sense that \eqref{eq:necon} holds.
\end{theorem}
We remark that the assumption on $R^TR$ is mild as $R^TR$ is indeed the estimated covariance matrix of assets. It is reasonable that a covariance matrix is positive definite when $N>d$, i.e., the number of realizations exceeds the number of assets. If we relax this assumption, we can add a proximal term in the $\x$-subproblem, which modifies our algorithm to a proximal ADMM algorithm like \cite{bai2021augmented}, and obtain a similar convergence result. For simplicity, we omit details for this case.
The assumption that $\cX$ is a convex and closed set is widely used  in practice, e.g., in many portfolio optimization applications, the weight has a box constraint  or a simplex constraint \citep{cui2018portfolio}.

As a final remark,  we point out that the decision making problem \eqref{pb:1} becomes easier to solve when $U$ is a strictly concave function or an identity function. This covers the cases of the RDU,  VaR, CVaR, SRM and DRM. Comparing to RDU, CVaR, SRM problems, VaR and DRM problem are much harder. However, our ADMM framework can handle these problems as well. Indeed, the $\x$-subproblems are still quadratic programs and can be solved by off-the-shelf solvers.
The $\y$-subproblems for above problems can be shown to be equivalent to convex separable problems with chain constraints and thus can be solved by either off-the-shelf solvers or specific algorithms for convex optimization with chain constraints.
For further details, please refer to \cref{rem:1} in the next section.

\section{Solving the $\y$-Subproblem in the ADMM}\label{sec:4}
In this section, we propose two methods to solve the $\y$-subproblem \eqref{eq:suby}. Before presenting our methods, we first introduce a reformulation for the $\y$-subproblem. Let $\w^{k+1} = R\x^{k+1} - {\blambda^k}/{\sigma}$. For simplicity, we omit the superscript in $\w^{k+1}$ in the following of this section. Let $\{l_1,l_2,\dots,l_N\}$ be a permutation of $\{1,\ldots,N\}$ such that
$
        w_{l_1} \leq w_{l_2} \leq \dots \leq w_{l_N}.
$
As $c_iU(y_{[i]})$ is independent on any permutation of $\y$,  the optimal solution $\y$ must be of the same ranking with $\w$, and hence the $\y$-subproblem \eqref{eq:suby} is equivalent to
\begin{equation}\label{eq:yisot}
        \min_{\y} ~\sum_{i=1}^N -c_iU(y_{l_i}) + \frac{\sigma}{2}(y_{l_i} - w_{l_i})^2  \quad
        \text{s.t.} \quad y_{l_1} \leq y_{l_2} \leq \dots\le y_{l_N}.
\end{equation}
The constraint $ y_{l_1} \leq y_{l_2} \leq \dots\le y_{l_N}$ is known as the simple chain constraint \citep{best2000minimizing}. Such a constraint first occurred in isotonic regression \citep{ayer1955empirical} and has been widely studied in the literature \citep{brunk1972statistical,stromberg1991algorithm}.
We remark that a similar reformulation idea was previously used in \cite[Lemma 3]{cui2018portfolio}.

\begin{remark}\label{rem:1}
When the function $U$ is a smooth and strictly concave function, or the identity function, which covers the cases of RDU, VaR, CVaR, SRM, and DRM,  problem \eqref{eq:yisot} becomes a convex optimization problem with a separable objective function over the chain constraint, which can be solved by a general convex solver. We can further utilize the celebrated pooling-adjacent-violators (PAV) method (e.g., \cite{best2000minimizing}) or the dynamic programming (DP) method proposed by \cite{yu2022dynamic}, which takes advantage of the structure of \eqref{eq:yisot}, to obtain much faster speed in both theory and practice. Therefore, we focus on the CPT utility function whose subproblem in \eqref{eq:yisot} is currently unsolvable by existing methods.
\end{remark}

Without loss of generality, in this section we assume $w_1 \leq w_2 \leq \dots \leq w_N$ for further simplicity of notations. Thus problem \eqref{eq:yisot} is equivalent to the following isotonic program
\begin{equation}\label{eq:subyiso}
        \min_{\y} ~\sum_{i=1}^N -c_iU(y_{i}) + \frac{\sigma}{2}(y_{i} - w_{i})^2  \quad         \text{s.t.}  \quad y_{1} \leq y_{2} \leq \dots \le y_{N}.
\end{equation}
Let $f_i(y_i) = -c_iU(y_i) + \frac{\sigma}{2}(y_i - w_i)^2$. That is,
$
f_i(y_i)=\left\{\begin{aligned} - a_i  U(y_i) + \frac{\sigma}{2}(y_i - w_i)^2  ,\quad  y_i \leq B,  \\ -b_i U(y_i) + \frac{\sigma}{2}(y_i - w_i)^2 , \quad y_i > B. \end{aligned}\right.
$
Note that according to $c_i$ in \eqref{eq:pd}, we have $a_i>0$ and $b_i>0$.

Then  problem \eqref{eq:subyiso} can be rewritten as
\begin{equation}\label{eq:dpyf}
        \min_{\y}~ \sum_{i=1}^N f_i(y_i) \quad
        \text{s.t.}  \quad y_1 \leq y_2 \leq \dots        \leq y_N.
\end{equation}
We will utilize the structure of $f_i$ to design our  algorithm.

To this end, we first consider the property of the following function
 \begin{equation}\label{eq:zeta}
\zeta(y)=\left\{\begin{aligned} & - a U(y)+ \frac{\sigma}{2}(y - w)^2,&  y \leq B,  \\
&-b U(y)+ \frac{\sigma}{2}(y - w)^2 , & y > B ,\end{aligned}\right.
\end{equation}
where $a,b \geq 0 ,\sigma >0$.
Note that  $f_i$  is exactly in the form of $\zeta$. It is obvious that $\zeta(y)$ is characterized by parameters $a,b,w,\sigma$, and function $U$. We will sometimes write it in the form $\zeta_{a,b,w,\sigma,U}(y)$ to emphasize the role of $a,b,w,\sigma$, and $U$.  The first order and  second order derivatives  of $\zeta(y) $ are given by
\begin{equation}\label{d1}
\zeta{'}(y)=\left\{\begin{aligned}&-  a U{'}(y)   + \sigma(y - w) , &y \leq B,  \\
&- b  U{'}(y) + \sigma(y - w) , & y > B ,\end{aligned}\right.
\quad\text{ and }\quad
\zeta{''}(y)=\left\{\begin{aligned} &- a  U{''}(y) + \sigma , & y \leq B,  \\
&- b U{''}(y) + \sigma ,&y > B .\end{aligned}\right.
\end{equation}

We next introduces properties of $\zeta$ that will be frequently used in the paper.
\begin{fact}\label{fact}
Assume $a>0,~b>0,$ and $\sigma>0$.

We first show that $\zeta(y)$ has at most one local minimizer and no local maximizer in interval $(B,\infty)$. Recall that $U(y)$ satisfies the properties in \cref{asmp:1}. The function $\zeta(y)$ is strongly convex in   $(B,\infty)$ due to $U{''}(y) < 0,~b>0$ and $\sigma > 0$, and $\zeta{'}(y)$ is strictly increasing in $(B,\infty)$. We have two cases depends on  the value of $\zeta{'}(B^{+}).$
\begin{enumerate}
	\item If $\zeta{'}(B^{+}) \geq 0$, then $\zeta{'}(y) > \zeta{'}(B^{+})\geq 0$ and thus $\zeta(y)$ is strictly increasing in $(B, \infty)$.
	\item If $\zeta{'}(B^{+}) < 0$, we have
	\begin{equation*}
		\zeta{'}(u) \ge -bU{'}(B+1) + \sigma( u - w) \geq -bU{'}(B+1) + \sigma \frac{bU{'}(B+1)}{\sigma} = 0,
	\end{equation*}
where $u = \max \left( B+1, w + \frac{bU{'}(B+1)}{\sigma}\right)$, and the first inequality follows from that $\zeta{'}(y)$ is strictly increasing in $(B,\infty)$.
As $\zeta{'}(B^{+}) < 0$ and $\zeta{'}(u) \ge 0$, there { exists} a unique $\xi \in (B,u]$ such that $\zeta{'}(\xi) = 0$. This, together with the fact that $\zeta'(y)$ is strictly increasing in $(B,\infty)$, yields that $\zeta(y)$ is strictly decreasing in $(B,\xi)$ and strictly increasing in $(\xi,\infty)$. We thus obtain that $\xi$ is a local minimizer.
\end{enumerate}

Next we show $\zeta(y)$ has at most one local minimizer and/or one local maximizer in interval $(-\infty, B)$. Note that $\zeta{'''}(y) = -a U{'''}(y) < 0$. This situation is more complicated.
\begin{enumerate}
	\item If $\zeta{''}(B^-) \geq 0$, then we have $\zeta{''}(y) > \zeta{''}(B^-) \geq 0, ~\forall y \in (-\infty,B)$. Hence $\zeta{'}(y)$ is strictly increasing in $(-\infty,B)$. Moreover,  we have $\zeta{'}(y) = -aU{'}(y) + \sigma(y-w) < 0$ for $y \le w$ as $a>0$ and $U'(y)>0$. Then we have the following two subcases.
   \begin{enumerate}
   	\item If $\zeta{'}(B^-) \leq 0$, then $\zeta(y)$ is strictly decreasing in $(-\infty,B)$.
   	\item If $\zeta{'}(B^-) > 0$, then we must have $B>w$, and there exists a unique $\xi \in (w, B)$ such that $\zeta{'}(\xi) = 0$. Furthermore, $\zeta(y)$ is strictly decreasing in $(-\infty,\xi)$ and increasing in $(\xi,B)$, and $\xi$ is a local minimizer.
   \end{enumerate}

   \item  If $\zeta{''}(B^-) < 0$, there exist two scenarios.
   \begin{enumerate}
   \item If $\lim_{y\rightarrow -\infty} \zeta{''}(y) > 0$, then there exists $C \in (-\infty,B)$ such that $\zeta{''}(C) = 0$,  $\zeta''(y)>0$ in $(-\infty ,C)$ and  $\zeta''(y)<0$  in $(C,B)$.
  We have the following two subcases.
  \begin{enumerate}
    \item $\zeta{'}(C) \leq 0$. It is easy to verify that $\zeta{'}(y) < 0,~ y \in (-\infty,C)\cup(C,B)$. Therefore $\zeta(y)$ is strictly decreasing in $(-\infty,B)$.
    \item $\zeta{'}(C) > 0$.  We must have $C>w$ as $U'(y)>0$ for $y<B$, $~a>0$ and $\sigma > 0$. There exists unique $\xi \in (w, C)$ such that $\zeta{'}(\xi) = 0$. One can verify that $\zeta(y)$ is strictly decreasing in $(-\infty,\xi)$ and increasing in $(\xi,C)$. So $\xi$ is a local minimizer. Next we consider $(C,B)$.
If $\zeta{'}(B^-) \geq 0$, then $\zeta{'}(y) > 0$ for $y \in (C,B)$ due to $\zeta''(y)<0$ for  $y\in(C,B)$, and thus $\zeta(y)$ is strictly increasing for $y \in (C,B)$. 
    Otherwise if $\zeta{'}(B^-)  < 0 $, there exists a unique $\pi \in (C,B)$ such that $\zeta{'}(\pi) = 0$ because $\zeta{''}(y)<0$ for all $y<B$ and $\zeta'(C)>0$, and thus $\pi$ is a local maximizer.
    \end{enumerate}
    \item  If $\lim_{y\rightarrow -\infty} \zeta{''}(y) \leq 0$, then we have $\zeta{''}(y) <0, \forall y \in (-\infty,B)$ as $\zeta'''(y)<0$ for $y\in(-\infty,B)$. We have $\lim_{y\rightarrow -\infty} \zeta{'}(y) = -a \lim_{y\rightarrow -\infty}U{'}(y) + \sigma(y-w) \le  0 + \sigma(y-w) < 0 $. Then we have $\zeta{'}(y) < 0,~\forall y \in (-\infty,B)$. Therefore, $\zeta$ is decreasing in $(-\infty,B)$.
\end{enumerate}

\end{enumerate}

It is possible that $B$ is either a local minimizer or local maximizer, depending on its monotonicity in $B$'s neighborhood.
\end{fact}

Note that finding a local minimizer (maximizer, respectively) of a univariate convex (concave, respectively) and continuously differentiable function in the interval $(\ell,u)$, which is assumed to exist, is equivalent to finding a zero point of its derivative. This is a root-finding problem for monotone and continuous functions, and can be easily solved by various methods like binary search or Fibonacci search to very high accuracy; see, e.g., \cite[Section 8.2]{bazaraa2006nonlinear}.
In our implementation, we adopt binary search that maintains a search interval including the optimal solution and perform one function evaluation to reduce the interval by half at each iteration.
We also remark that assuming $G(\theta)$ is a monotone and continuous function in interval $(\ell,u)$, it takes at most $O\left(\log_2(\frac{u-\ell}{\epsilon})\right)$ function evaluations to find
an approximate root $\bar\theta$ such that $|\bar\theta-\theta^*|<\epsilon$, where $\theta^*$ is a root for $G(\theta)=0$ in $[\ell,u]$ \citep{ahuja2001fast}. In practice, several tens of binary searches are sufficient to find a solution up to machine accuracy.
Therefore, in this paper we always assume that there is an oracle that finds the exact root of a monotone and continuous function.
\begin{definition}[Root-finding oracle]
A root-finding oracle is a procedure that either finds the root or claims that there is no root for a univariate function that is monotone and continuous in a given interval.
\end{definition}

Therefore, if $F(\theta)$ is a convex (concave, respectively)  continuously differentiable function in $(\ell,u)$, then the root-finding oracle for $F'(\theta)$ can find a local minimizer (maximizer, respectively).
Using this and \cref{fact}, we conclude the following lemma.
\begin{lemma}\label{lem:3ora}
Consider $\zeta$ in \eqref{eq:zeta} with $a,b,\sigma>0$ and $U$ satisfying \cref{asmp:1}.
It takes at most three root-finding oracles to find all the local minimum for $\zeta$ in an arbitrary interval in $(-\infty,\infty)$.
\end{lemma}

We next present a useful result that the optimal solution of the $\y$-subproblem must satisfy a condition so that the Clarke generalized gradient is well defined according to \cref{lem:locallip}, which plays an essential role in our algorithm design and the associated convergence analysis.
\begin{lemma}\label{lem:4}
Let $\y^*$ be an optimal solution of problem \eqref{eq:dpyf}. Then we have $U'(y^*_i)< \infty,~\forall i = 1,2,\dots,N$.
\end{lemma}

The next lemma shows that the optimal solution of problem \eqref{eq:dpyf} is  bounded by known constants. This is a useful property for the convergence analysis of our algorithm.
\begin{lemma} \label{lem:bound}
Let $\y^*$ be an optimal solution of problem\eqref{eq:dpyf}. Then we have
\[
w_1\le y_1^*\le \dots \le y_N^* \le \max_{i=1,\ldots,N}\left\{  \max \left( B+1, w_i + \frac{b_iU'(B+1)}{\sigma}\right)\right\}.
\]
\end{lemma}
In the following of this paper, we set
\[l_b=w_1\quad\text{and}\quad u_b=\max_{i=1,\ldots,N}\left\{  \max \left( B+1, w_i + \frac{b_iU'(B+1)}{\sigma}\right)\right\}.\]
The above $l_b$ and $u_b$ serve as initial lower and upper bounds for the solutions for the following proposed algorithms for the $\y$-subproblem.

\subsection{Dynamic Programming Algorithm}\label{sec:dp}
This section employs the philosophy of dynamic programming (DP) to present an algorithm that provides a global solution to problem \eqref{eq:dpyf}. Our approach is inspired by the DP technique for generalized nearly isotonic optimization proposed in \cite{yu2022dynamic}. However, it should be noted that the method in \cite{yu2022dynamic} is exclusively appropriate for convex separable objective functions and is thus not applicable to our situation. The DP method's principal advantage is that it yields an exact solution for our nonconvex objective function by comprehensively leveraging the problem's structure.

To introduce our algorithm, let us first define $h_1(z) = 0$, and define $h_n(z)$ recursively as
\begin{equation}\label{eq:hupdate}
        h_n(z) = \min_{y \leq z} f_{n-1}(y) + h_{n-1}(y), \quad n=2,\ldots,N+1.
\end{equation}
By definition, $h_n(z)$ is a non-increasing function of $z$ as a bigger $z$ enlarges the feasible region of the minimization problem in \eqref{eq:hupdate}.
Let $\y^*$ be an optimal solution of \eqref{eq:dpyf}.
Noting that \eqref{eq:hupdate} is equivalent to
\[
h_{N+1}(y_{N+1})=    \min_{y_N\le y_{N+1}}\left\{ f_N(y_N)+ \min_{y_{N-1}\le y_N} \left \{ f_{N-1}(y_{N-1})+\cdots +\min_{y_1\le y_2}\left\{f_1(y_1)+h_1(y_1)\right\}\right\}\right\},
\]
we can verify that \eqref{eq:dpyf} is equivalent to
\begin{align}
\min_{u\le y_{N+1}}~ f_N(u) + h_N(u)  \label{finalobj}
\end{align}
for any given $y_{N+1}\ge y_N^*$.
Problem \eqref{finalobj} can be seen as a univariate unconstraint optimization problem.

An outline of our DP algorithm is that we first forwardly derive the explicit expressions for $h_2,h_3,\ldots,h_{N+1}$, and then use the expression
\begin{equation}\label{eq:hrecur}
y_{n-1}^* = \mathop{\text{argmin}}_{y\leq y_n^*}f_{n-1}(y) + h_{n-1}(y),\quad n=N+1,N,\ldots,2
\end{equation}
 backwardly to recover the optimal $y_N^*,y_{N-1}^*,\ldots,y_1^*$.

\subsubsection{Forward Update of the DP Algorithm}
A key observation of the DP algorithm is that $h_n$ is a piecewise function that in each piece, $h_n$ is either a constant, or a strictly decreasing function in the form of $\zeta(y)$ plus a constant.
Based on this, the DP algorithm applies a forward update to derive the expression of $h_n$, using the expression of $h_{n-1}$ and decomposing the interval $[l_b,u_b]$ into consecutive intervals so that the expression (i.e. the parameters) of $h_{n}$ is unified in each interval for $n=1,2,\ldots,N+1$.
In the following, we will show that $h_{n}(z)$ is a piecewise function with $K_{n}$ consecutive intervals, say $[s^{n}_k,s^{n}_{k+1}],~k = 1,2,...,K_{n}$, and these intervals form $[l_b,u_b]$.
We will also obtain the explicit expression of $h_{n}$ in each piece in a unified form. Let $L_n$ denote the set containing these intervals.

To find the expression of $h_n(z) $ in  $[s^{n-1}_k,s^{n-1}_{k+1}]$, we need first explore the explicit formula for an auxiliary function  defined by
\begin{equation}\label{eq:phiy}
        \phi_{n-1}(z) = f_{n-1}(z) + h_{n-1}(z) ,\quad z\in  [s^{n-1}_k,s^{n-1}_{k+1}], ~ k = 1,2,...,K_{n-1}.
\end{equation}
The above formula, together with \eqref{eq:hupdate}, yields
\begin{equation}\label{eq:hleqphi}
        h_n(z) = \min_{x \leq z} \phi_{n-1}(x), \quad z\in [s^{n-1}_k,s^{n-1}_{k+1}],~k=1,2,\ldots,K_{n-1}.
\end{equation}
Once we know the expression of $h_{n-1}(z)$ in intervals $[s^{n-1}_1,s^{n-1}_{2}],\ldots,[s^{n-1}_{k-1},s^{n-1}_{k}]$,
we can derive the expression of $h_n(z)$ in $[s^{n-1}_k,s^{n-1}_{k+1}]$ by
\begin{align}
        h_n(z) &= \min \left(\min_{x\leq s^{n-1}_{k}} \phi_{n-1}(x),  \min_{s^{n-1}_{k}\leq x \leq z} \phi_{n-1}(x)\right) , \quad z \in [s^{n-1}_k,s^{n-1}_{k+1}], \notag\\
                &= \min\left(\ h_n(s^{n-1}_{k}), \min_{s^{n-1}_{k}\leq x \leq z} \phi_{n-1}(x)\right), \quad z \in [s^{n-1}_k,s^{n-1}_{k+1}],\label{hphi}
\end{align}
where the second equation follows directly from the definition of $h_{n-1}$ in \eqref{eq:hupdate}.
Repeating the update \eqref{hphi} for each interval $[s^{n-1}_k,s^{n-1}_{k+1}]$, we obtain the full expression for $h_n$ in the interval $[l_b,u_b]$.
Next, based on \eqref{eq:hleqphi} we will show that $h_n(z)$ is a piecewise function,
which is either strictly decreasing or constant in each piece.

To facilitate our analysis, we give a lemma on the continuity of $h_n(z)$, which follows from \eqref{eq:hupdate}, the continuity of $f_i$, and $h_1(z)=0$.
\begin{lemma}\label{hncont}
The function $h_n(z)$ is a continuous function in $(-\infty,\infty)$.
\end{lemma}

\begin{center}
\begin{minipage}{\linewidth}
   \begin{algorithm}[H]\small
        \caption{The deduction of $h_n(\cdot)$ from $h_{n-1}(\cdot)$ in $[l_b, u_b]$}
        \label{dpalg1}
        \begin{algorithmic}[1]
            \Require$h_{n-1}(\cdot), f_{n-1}(\cdot)$, and $ L^{n-1}=\{ [s^{n-1}_{1},s^{n-1}_{2}],[s^{n-1}_{2},s^{n-1}_{3}],\dots, [s^{n-1}_{K_{n-1}^{L}},s^{n-1}_{K_{n-1}^{L}+1}]\}$  with $s^{n-1}_{1}=l_b$ and $s^{n-1}_{K_{n-1}+1}=u_b$
            \Ensure $\{ h_n(\cdot), L^n \}$
            \State $h_n(s^{n-1}_{1}) = \phi_{n-1}(s^{n-1}_{1})$; $L^{n} = \emptyset$\Comment{$\phi$ is defined in \eqref{eq:phiy}}

            \For {$k=1,2,\dots,K_{n-1}$}
            \State $T^{n-1}_{k} = \Call{Decompose}{[s^{n-1}_k,s^{n-1}_{k+1}], \phi_{n-1}(\cdot)}$
                \For {$i=1,\ldots,\Card(T^{n-1}_{k})$}
\Comment{$\Card(T^{n-1,L}_{k})$ is the cardinality of $T^{n-1}_{k}$}
\State Select $[t_i, t_{i+1}] \in T^{n-1}_{k}$;
              \State $\{Q_i, h_n(\cdot)\} = \Call{Update}{[t_i, t_{i+1}], \phi_{n-1}(\cdot), h_{n}(t_i)}$;
               \State $L^{n} = \Call{Merge}{L^{n}\cup Q_i}$
                \EndFor
            \EndFor
        \end{algorithmic}
\end{algorithm}
\end{minipage}
\end{center}

To better explore the expression of $h_n(z)$, we introduce three key procedures. Provided $h_{n-1}$, the update of $h_n$ is summarised in \cref{dpalg1}. The first procedure is used to decompose each piece of $\phi_{n-1}=f_{n-1}+h_{n-1}$ to smaller consecutive intervals with desirable properties.
\begin{procedure}[$\Call{Decompose}{}$]\label{def:decom}
For a function in the form of $\zeta(\cdot)$ with $a,b,\sigma>0$ in an interval $[p,q]\subset [l_b,u_b]$. If $\zeta(\cdot)$ is neither strictly decreasing nor strictly increasing in $[p,q]$, then we decompose $[p,q]$ into several smaller intervals where in each new interval $\zeta(\cdot)$ is either strictly increasing or strictly decreasing. We refer to this procedure as $\Call{Decompose}{}$.
It outputs a set of consecutive intervals, i.e., $T =\Call{Decompose}{[p,q], \zeta(\cdot)}$.
\end{procedure}
The well-definedness of \Call{Decompose}{} follows from the structure of $\zeta(\cdot)$.
Now let us illustrate how $\Call{Decompose}{}$ works. Suppose $\zeta(y)$ is  neither strictly decreasing nor strictly increasing in $[p,q]$. From \cref{fact}, we can find the local minimizer and/or maximizer of $\zeta(y)$ in $[p,q]$, using which as end points, the interval $[p,q]$ can be decomposed into at most three consecutive intervals. In each of these intervals, $\zeta(y)$ is either strictly decreasing or increasing.
We remark that the local minimizer or maximizer corresponds to the root of $\zeta'(\cdot)$ in the associated interval, and thus can be found by our root-finding oracle.

The next procedure updates the expression of $h_n$ using the philosophy of \eqref{hphi}.
\begin{procedure}
[$\Call{Update}{}$]
Let $\delta(\cdot)$ be a strictly monotone and continuous function defined in the interval $[p,q]$ and $\kappa$ be a constant satisfying $\kappa\leq \delta(p)$. Define $\theta(z) = \min\{ \kappa, \min_{p \leq y \leq z} \delta(y) \}$. The procedure $\Call{Update}{}$ returns an explicit piecewise expression of $\theta(z)$ in the interval $[p,q]$. It works as $\{S, \theta(\cdot)\} = \Call{Update}{[p,q], \delta(\cdot), \kappa}$, where $S$ denotes the consecutive intervals that support $\theta(\cdot)$.
\end{procedure}
Now we describe the details of $\Call{Update}{}$.
When $\delta(y)$ is strictly increasing, because of $\min_{p \leq y \leq z} \delta(y) = \delta(p) \geq \kappa$, from the definition of $\theta(z)$, we have $\theta(z) = \kappa, z\in [p,q]$, and
 $S = \{ [p,q] \}$. The case that  $\delta(y)$ is a constant is trivial.
We next show all possible outputs of this procedure in the nontrivial case that $\delta(y)$ is strictly decreasing. There are three scenarios.
\begin{enumerate}
        \item If $\delta(p) = \kappa$, then  $\theta(z) = \delta(z),~ z\in [p,q]$ and  $S = \{ [p,q] \}$.
        \item If $\delta(q) \geq \kappa$, then   $\theta(z) = \kappa, ~z\in [p,q]$,  and $S = \{ [p,q] \}$.
    \item If $\delta(p) > \kappa > \delta(q)$, then $\theta(z)$ is a piecewise function with two pieces. Indeed, the strictly decreasing property and continuity of $\delta(\cdot)$ guarantees that exists a unique point in $(p,q)$, denoted by $r$, such that $\delta(r) = \kappa$.

Thus we have 
$\small
\theta(z)=\left\{\begin{array}{ll} \kappa ,&  z \in [p,r],  \\
 \delta(z) , & z \in [r,q], \end{array}\right.
$
and $S = \{ [p,r] ,[r,q]\}$.
\end{enumerate}
In Algorithm \ref{dpalg1},  \Call{Update}{}  is used to update the expression of $h_n$ (with $\delta=\phi_{n-1}$) in each piece. After \Call{Update}{}, $h_n(z)$ is a piecewise function,
which is either strictly decreasing or constant in each piece.
We also remark that   the time cost of  \Call{Update}{} is dominated by computing a root for a monotone equation, which can be done by our root-finding oracle.

Now let us show in the next lemma that $h_n$ admits a unified expression in each piece, by using the two procedures \Call{Decompose}{} and  \Call{Update}{}. The proof is deferred to Appendix \ref{proof:l7}.
\begin{lemma}\label{lem:recur}
For each $m=1,2,\ldots,N+1$ and $z\in[l_b,u_b]$, $h_{m}(z)$ can be characterized by a piecewise function with $K_{m}$ pieces, given by
        \begin{equation}\label{eq:hexpre}
        h_{m}(z) = \zeta_{a_k^{m},b_k^{m},w_k^{m},\sigma_k^{m},U}(z) + M_k^{m},~z\in [s^{m}_k,s^{m}_{k+1}], \quad ~k = 1,2,\dots,K_{m},
        \end{equation}
where $s^m_1=l_b$, $s^{m}_{K_{m}}=u_b$, $M_k^{m}$ is constant, and either $a_k^{m} > 0,~b_k^{m}>0,~\sigma_k^{m} > 0$ or $a_k^{m} = b_k^{m} =\sigma_k^{m}=0$.
\end{lemma}

The final main issue is that once  $h_n(y)$ is constant in two consecutive  intervals after the \Call{Update}{} procedure, we need to merge the two intervals so that the following backtracking procedure, which will be used to recover the optimal solution $\y^*$, is well defined.

\begin{procedure}[\Call{Merge}{}]
Suppose that after the \Call{Update}{} procedure,  $h_n(y_n)$ is  constant in two consecutive intervals $[s^n_{i-1},s^n_{i}]$ and $[s^n_{i},s^n_{i+1}]$. Note that due to continuity, the constant values of $h_n$ in the two intervals are the same. Then we merge the two intervals into $[s^n_{i-1}, s^n_{i+1}]$ and obtain that $h_n$ is constant in $[s^n_{i-1}, s^n_{i+1}]$ .
The \Call{Merge}{} procedure merges all consecutive intervals where $h_n$ is constant.
It works as $S=\Call{Merge}{T}$, where $T$ consists of consecutive intervals.
\end{procedure}
We remark that the case that $h_n$ is not constant but admits the same expression in two consecutive pieces cannot occur. Indeed, by our update rule,  in each piece $h_n$ is determined by the constant $M_k^{m}$ in $h_k$ and $\sum_{i=k}^{n-1}f_i(\cdot)$ for some $k \leq n-1$.
If the two consecutive pieces share the same $\sum_{i=k}^{n-1}f_i(\cdot)$, then the constants must be the same because of the continuity of $h_n$. This contradicts the \Call{Merge}{} step as we have merged the consecutive constant intervals in $h_{k}$ for $k<n$.

\subsubsection{Backward Procedure of the DP Algorithm}
In the last subsection,  we obtained a forward procedure to derive the explicit formula for $h_n(\cdot)$, $n = 1,2,\dots, N+1$. We still need a backward procedure to obtain optimal values for $y^*_N, y^*_{N-1},\dots, y^*_2,y^*_1$. To this end, we introduce several lemmas below, whose proof follows directly from \eqref{eq:hleqphi}.

\begin{lemma}\label{hndecrease}
If  $h_n(z)$ is strictly decreasing in $[s^{n}_{k},s^{n}_{k+1} ]$,  then $h_n(z) = \phi_{n-1}(z)$ for $z\in [s^{n}_{k},s^{n}_{k+1} ]$.
\end{lemma}

\begin{lemma}\label{lem:hnconstant}
Consider the stage after \Call{Merge}{}. Suppose that $h_n(z)$ is a constant in $[s^{n}_{k},s^{n}_{k+1}]$. Then $h_n(z) = \phi_{n-1}(s^{n}_{k})$ for all $z\in[s^{n}_{k},s^{n}_{k+1}]$, and thus $\phi_{n-1}(y) \geq \phi_{n-1}(s^{n}_{k})$ for all $y\in [s^{n}_{k},s^{n}_{k+1}]$.
 \end{lemma}

Now we are ready to present a backward recursion for obtaining the optimal solution $\y^*$.
The original problem  \eqref{eq:dpyf} is equivalent to \eqref{finalobj}.
Using \eqref{eq:hrecur}, we can recursively update $y_n^*$ backward,
 by setting $y_{N+1}^*=u_b $. Now we will show the details of  how to find $y_{n-1}^*$ using $y_{n}^*$ for $n=N+1,N,\ldots,2$.
Suppose $[s^{n}_{k},s^{n}_{k+1}]$ is the piece of $h_n$ that contains $y_n^*$. We update $y_{n-1}^*$ according to the following two cases.
\begin{enumerate}
\item  $h_n(y)$ is strictly decreasing in the interval $[s^{n}_{k},s^{n}_{k+1}]$.
Then $h_n(y) =\phi_{n-1}(y), ~y \in [s^{n}_{k},s^{n}_{k+1}]$ due to Lemma \ref{hndecrease}. So $\phi_{n-1}(y)$ is strictly decreasing in  $[s^{n}_{k},s^{n}_{k+1}]$. Particularly, $\phi_{n-1}(s^{n}_{k}) = h_n(s^{n}_{k})$. Therefore $h_n(s^{n}_{k}) \geq \phi_{n-1}(y), \forall y \in [s^{n}_{k}, y^*_n]$. From \eqref{eq:hleqphi}, we have
\[ h_n(s^{n}_{k}) =\min_{x\le s^{n}_{k}} \phi_{n-1}(x)\leq \phi_{n-1}(z) , \quad \forall z \leq s^{n}_{k}.\]
The above two facts imply $\phi_{n-1}(z)\ge \phi_{n-1}(y)$ for all $ z \leq s^{n}_{k}$ and $y\in [s^{n}_{k}, y^*_n]$
and thus
\begin{equation*}
        y_{n-1}^* \overset{\eqref{eq:hrecur},~\eqref{eq:phiy}}{=} \mathop{\text{argmin}}_{y\leq y_n^*} \phi_{n-1}(y) = \mathop{\text{argmin}} _{s^{n}_{k}\leq y\leq y_n^*}\phi_{n-1}(y) = y_n^*.
\end{equation*}
The last equation is because $\phi_{n-1}(y)$ is strictly decreasing in $[s^{n}_{k},y^*_n]$.

\item  $h_n(y)$ is constant in the interval $[s^{n}_{k}, s^{n}_{k+1}]$. This means $\phi_{n-1}(y)\geq h_n(s^{n}_{k}) = \phi_{n-1}(s^{n}_{k}), y \in [s^{n}_{k},y_n^*]$, thanks to Lemma \ref{lem:hnconstant}. This, together with \eqref{eq:hleqphi}, guarantees that $\phi_{n-1}(s^{n}_{k}) = h_n(s^{n}_{k}) \leq \phi_{n-1}(y),~ \forall y \leq s^{n}_{k}$. Therefore as $h_n(y)$ is a constant in the interval $[s^{n}_{k}, y_n^*]$, we have
\begin{equation*}
        \phi_{n-1}(s^{n}_{k}) = h_n(s^{n}_{k}) = h_n(y_n^*)= \min_{y_{n-1}\leq y_n^*} \phi_{n-1}(y_{n-1}).
\end{equation*}
This implies
$
        y_{n-1}^* = \mathop{\text{argmin}}_{y_{n-1}\leq y_n^*} \phi_{n-1}(y_{n-1}) = s^{n}_{k}.
$
\end{enumerate}
\begin{center}
\begin{minipage}{\linewidth}
\begin{algorithm}[H]\small
        \caption{A dynamic programming algorithm for solving \eqref{eq:dpyf}}
        \label{alg:dp}
        \begin{algorithmic}[1]
            \Require $f_i$, $i=1,\ldots,N$, $h_1(z) = 0$, $L^1 = \{[l_b,u_b]\}$
            \For {$n = 2,3,\dots, N+1$}
            \State Obtain $\{ h_n(\cdot), L^n \}$ by  \cref{dpalg1}
        \EndFor
     \State $y_{N+1}^* = u_b$
     \For{$n = N+1,N,\dots,2$}
     \State Find $[s^{n}_{k},s^{n}_{k+1}] \in L_n$ such that $y_n^* \in [s^{n}_{k},s^{n}_{k+1}]$
     \If{$h_n(\cdot)$ is a constant in $[s^{n}_{k},s^{n}_{k+1}]$}
     \State $y_{n-1}^* = s^{n}_{k}$
     \Else
     \State $y_{n-1}^* = y_n^*$
     \EndIf
     \EndFor
        \end{algorithmic}
\end{algorithm}
\end{minipage}
\end{center}

\subsubsection{The DP algorithm}
We summarize the full DP algorithm in \cref{alg:dp}.
We remark that in line 3, we introduce the dummy variable $y_{N+1}^* = u_b$ according to \eqref{finalobj}.
The following theorem is then straightforward.
\begin{theorem}\label{thm:dp}
\cref{alg:dp} correctly returns a global optimal solution for problem \eqref{eq:subyiso}.
\end{theorem}

The time complexity of our DP algorithm mainly depends on the number of pieces of intervals in each $h_n$. Assume this number is $O(M)$ for every $n=1,2,\dots,N$. The intervals are generated by the $\Call{Decompose}{}$ and $\Call{Update}{}$ steps. It is worth recalling that each $\Call{Decompose}{}$ or $\Call{Update}{}$ step triggers the root-finding oracle at most three times, as affirmed by \cref{lem:3ora}. Considering the existence of $N$ stages in the forward deduction, the DP algorithm executes the oracle at most $O(MN)$ times.  On the other hand, in the backward recursion, at most $O(MN)$  intervals need to be traversed, which is faster than the forward stage as no root-finding oracle needs to be invoked. Therefore the DP algorithm invokes the oracle  $O(MN)$ times.
In the worst-case scenario, an interval may generate up to five new consecutive intervals in each iteration after the $\Call{Decompose}{}$ and $\Call{Update}{}$ steps. This is attributable to the presence of at most one local minimizer in $(B,\infty)$, one local minimizer and one local maximizer in $(-\infty,B)$, and the possibility of the point $B$ constituting an extreme point. This provides an upper bound for $M$, i.e., $M \leq 5^N$.
However, in our (thousands of) numerical experiments for the y-subproblem in CPT utility maximization problem, $M$ typically ranges around $N$. Our empirical observations also suggest that $h_{N+1}$ comprises the most intervals among $h_i,~i=1,\ldots,N+1$. 
In practice, we may roughly think that the oracle complexity of our DP method is $O(N^2)$.

\subsection{Pooling-adjacent-violators Algorithm}
The DP method returns a global optimal solution for the $\y$-subproblem, but it suffers from a relatively low speed. To remedy this, we propose a much faster algorithm for the $\y$-subproblem, which computes a stationary point of a reformulation of the $\y$-subproblem \eqref{eq:suby}.

The pooling-adjacent-violators (PAV) algorithm is widely regarded as one of the most successful methods for resolving separable convex optimization subject to the chain constraint. Initially, the PAV algorithm was proposed for maximum likelihood estimation (MLE) under the chain constraint in \cite{ayer1955empirical} and \cite{brunk1972statistical}.

 The PAV methods proposed in \cite{stromberg1991algorithm,best2000minimizing} and \cite{ahuja2001fast} all require  $\theta_i(\cdot)$ to be a convex function, a condition that is not met in our problem. While the PAV-like method introduced in \cite{cui2021non} is tailored for a specific type of non-convex $\theta_i(\cdot)$, it imposes an additional requirement that $\theta_i(y) = k_i a(y) + m_i b(y) + l_i$, where $a(y)$ is strictly increasing and $m_i > 0$. Unfortunately, our objective function cannot be transformed into such a form.

Nonetheless, we present an algorithm founded on the PAV algorithm by exploiting the distinctive structure of $f_i(y)$, despite the nonconvexity of both \eqref{eq:suby} and the reformulation \eqref{eq:dpyf}.

Now let us first recall the traditional PAV algorithm, which aims to solve the following problem
\[
\min_{\y}~\sum_{i=1}^n \theta_i(y_i) ~ \text{ s.t. } \ y_1\leq y_{ 2} \leq \dots \leq y_n,
\]
where each $\theta_i$ is a univariate convex function.
The PAV algorithm maintains a set $J$ that partitions indices $\{1,2,\dots,n\}$ into consecutive blocks $\{[s_1+1,s_2],[s_2+1,s_3],\cdots,[s_k+1,s_{k+1}]\}$ with $s_1=0$ and $s_{k+1}=N$, where $[a,b]$ denotes the index set $\{a,a+1,\ldots,b\}$ for positive integers $a<b$, and we define by convention $[a,a]=\{a\}$.  A block $[p,q]$ is a \textit{single-valued block} if every $y_i$ in the optimal solution of the following problem has the same value,
\begin{equation*}
        \min_{y_p,\dots, y_q} ~ \sum_{i=p}^q \theta_i(y_i) \quad \text{s.t.} \ y_p\leq y_{p+1} \leq \dots \leq y_q,
\end{equation*}
i.e., $y_p^* = y_{p+1}^* = \dots = y_q^*$. We use $V_{[p,q]}$  to denote this value, following the literature.
For two consecutive blocks $[p,q]$ and $[q+1,r]$, if $V_{[p,q]} \leq V_{[q+1,r]}$ then the two blocks are \textit{in-order}, otherwise \textit{out-of-order}. The PAV algorithm initially partitions every  integer from 1 to $n$ as single-valued blocks $[i,i]$, $i=1,\ldots,n$. Once there exists consecutive out-of-order single-valued blocks $[p,q]$ and $[q+1,r]$, the PAV algorithm merges these two blocks by replacing them with the larger block $[p,r]$. When all the single-valued blocks are in-order, the PAV algorithm terminates.

Now let us describe our variant of the PAV algorithm.  Using the convention \[g_{[p,q]}(y)= \sum _{i = p}^q f_i(y),\] $g_{[p,q]}(y)$ can be parameterized as $\zeta(y)$ plus a constant by a similar discussion in the proof of \cref{lem:recur}. Despite that $g_{[p,q]}(y)$ is nonconvex, we can find its minimizer according to \cref{fact}. In fact, $g_{[p,q]}(y)$ has at most two local minimizers, we can simply choose the minimizer with the smaller function value. Our variant of the PAV is summarised in \cref{alg:pav}. However, the guarantee to a global optimal solution  of the PAV never holds now. Nevertheless, our numerical result suggests that the proposed algorithm has a good practical performance. 
\begin{center}
\begin{minipage}{\linewidth}
\begin{algorithm}[H]\small
        \caption{A  Pool-adjacent-violators algorithm for solving \eqref{eq:dpyf}}
        \label{alg:pav}
        \begin{algorithmic}[1]
            \Require  $J = \{ [1,1],[2,2],\dots,[N,N]\}$
            \For {each $[i,i] \in J$}
            \State Compute  the  minimizer $V_{[i,i]}$ of $f_i(y)$
        \EndFor
        \While {$\exists [m,n],[n+1,p] \in J$ such that $V_{[m,n]}>V_{[n+1,p]}$}
        \State  $J\leftarrow J\backslash\{[m,n],[n+1,p]\}\cup \{[m,p]\}$;       $g_{[m,p]}(y) \leftarrow \sum_{i=m}^pf_i(y)$;
    Compute  the  minimizer $V_{[m,p]}$ of $g_{[m,p]}(y)$
        \EndWhile
        \For {each $[m,n] \in J$ }
        \State $y_i = V_{[m,n]}, ~\forall i = m,m+1,\dots,n$
        \EndFor
           \end{algorithmic}
\end{algorithm}
\end{minipage}
\end{center}

Although we cannot expect the proposed variant of the PAV converges to a global optimal solution, we show in the next theorem that the PAV algorithm generates a stationary point for the following problem
\begin{equation*}\label{eq:eqsub}
        \min_{\y} ~\Gamma(\y) :=\quad \sum_{i=1}^N -c_iU(y_{[i]}) + \frac{\sigma}{2}(y_{[i]} - w_{i})^2,
\end{equation*}
which is equivalent to \eqref{eq:subyiso}.
\begin{theorem}
\label{thm:3}
        The solution $\bar\y$ returned by the PAV algorithm satisfies
$\bz\in \partial \Gamma(\bar \y).$
\end{theorem}
 We also remark that in our empirical test the descent property \eqref{eq:ydescent} is always satisfied.

Finally, we present the oracle complexity of the PAV algorithm.
\begin{theorem}\label{thm:4}
Given the root-finding oracle for computing minimizers in lines 2  and 4, \cref{alg:pav}  invokes the oracle $6N-3$ times.
\end{theorem}

We finally give a remark on the runtime of a function evaluation of each oracle in line 4.
In each merging step, the time complexity of calculating $g_{[p,r]}$ is $O(1)$ due to the special structure of $g_{[p,r]}$.
Indeed, in the merging step, deriving the coefficients of $g_{[p,r]}$ from $g_{[p,q]}$ and $g_{[q+1,r]}$ can be done in $O(1)$ time.
To see this, let $g_{[p,q]}(y) = \zeta_{a_{[p,q]},b_{[p,q]}, w_{[p,q]},\sigma_{[p,q]},U}(y)+ {\rm const}$ and $g_{[q+1,r]}(y) = \zeta_{a_{[q+1,r]},b_{[q+1,r]}, w_{[q+1,r]},\sigma_{[q+1,r]},U}(y)+ {\rm const}$, where ${\rm const}$ denotes the constant term. Then the parameters for $g_{[p,r]}$ are calculated as follows, $
        a_{[p,r]} = a_{[p,q]} + a_{[q+1,r]}$,
        $b_{[p,r]} = b_{[p,q]} + b_{[q+1,r]}$,
        $\sigma_{[p,r]} = \sigma_{[p,q]} + \sigma_{[q+1,r]}$,
        $w_{[p,r]}  = \frac{\sigma_{[p,q]]}w_{[p,q]}+ \sigma_{[q+1,r]} w_{[q+1,r]}}{\sigma_{[p,q]} + \sigma_{[q+1,r]}}$. Once the paramters in $g_{[p,r]}$ are calculated, the root-finding oracle can be used to find a local minimizer of $g_{[p,r]}$.
So we conclude the PAV algorithm takes $O(N)$ time from \cref{thm:4}.
We should point out that the existing PAV algorithm in \cite{ahuja2001fast,stromberg1991algorithm,best2000minimizing} requires $O(N)$ time to evaluate a general block objective function as it involves $O(N)$ additions of $f_i$, which is more expensive than $O(1)$ in our algorithm.

\section{Numerical Experiments}\label{sec:5}
This section presents the results of our performance testing of the proposed methods. Firstly, we compare the effectiveness of our two $\y$-subproblem solvers, namely the DP and PAV algorithms, on random instances. Next, we evaluate the performance of our ADMM methods on a portfolio optimization problem, with both the DP and PAV algorithms used as solvers for the $\y$-subproblem. Lastly, we conduct an empirical study for portfolio optimization using our algorithm.  All the experiments are implemented using MATLAB {R2023b} on a PC running Windows 10 Intel(R) Xeon(R) E5-2650 v4 CPU (2.2GHz) and 64GB RAM. The implementation of our code can be found in \href{https://github.com/RufengXiao/ADMM_CPT}{https://github.com/RufengXiao/ADMM\_CPT}.

\subsection{Numerical Tests for the $\y$-subproblem Solvers}\label{sec:5.1}
In this subsection, we aim to test the numerical efficiency of our proposed two $\y$-subproblem solvers when the utility function $U$ is given by \eqref{eq:sshapeu}. The parameters for the utility are $\mu=2.25, \alpha=0.88,\delta=0.69,\gamma=0.61$ and $B = 0$.
We compare the average performance of the DP and PAV algorithms with the \texttt{fmincon} solver in MATLAB; To prevent premature termination and ensure a more reliable comparison of objective function values, we configure the parameter \textsf{MaxFunctionEvaluations} to 1,000,000 and \textsf{MaxIterations} to 3,000. All other parameters are set to their default values.
We test different instances with  varying numbers of scenarios. We generate 10 instances for each scenario number, where each component of $\w$ in \eqref{eq:subyiso} is uniformly  drawn from $[-0.1,0.1]$.

\begin{center}
\begin{minipage}{\linewidth}
\begin{table}[H]  \caption{Objective value and time comparison for the $\y$-subproblem \eqref{eq:subyiso}.}
\label{table:objsub}
{\tiny
\begin{tabular*}{\textwidth}{l@{\extracolsep{\fill}}|lll|lll}
\toprule
\multirow{2}{*}{Scenarios } &\multicolumn{3}{c|}{objective value ($\times 10^{-2}$)}&\multicolumn{3}{c}{time (seconds)}\\
  &DP  & PAV & \texttt{fmincon}  & DP & PAV & \texttt{fmincon} \\
\midrule
$N=50$&\textbf{4.3058} $\pm$ 1.0884&\textbf{4.3058} $\pm$ 1.0884 &4.3127 $\pm$ 1.0807&0.0866 $\pm$ 0.0500&\textbf{0.0177} $\pm$ 0.0332 &0.5851 $\pm$ 0.7271\\
$N=100$&\textbf{4.6096} $\pm$ 0.7960&4.6097 $\pm$ 0.7959&4.6101 $\pm$ 0.7961&0.2358 $\pm$ 0.0296&\textbf{0.0095} $\pm$ 0.0021 &1.6821 $\pm$ 0.4736\\
$N=200$&\textbf{4.8137} $\pm$ 0.7721&\textbf{4.8137} $\pm$ 0.7721&4.8228 $\pm$ 0.7781&0.7757 $\pm$ 0.0504&\textbf{0.0180} $\pm$ 0.0034 &3.4323 $\pm$ 0.9501\\
$N=300$&\textbf{4.5241} $\pm$ 0.4999&4.5242 $\pm$ 0.4999&4.5386 $\pm$ 0.4968&1.6700 $\pm$ 0.0413&\textbf{0.0259} $\pm$ 0.0056 &6.6714 $\pm$ 1.7944\\
$N=400$&\textbf{4.8163} $\pm$ 0.5271&\textbf{4.8163} $\pm$ 0.5271&4.8191 $\pm$ 0.5243&2.9349 $\pm$ 0.0863&\textbf{0.0373} $\pm$ 0.0086 &9.9095 $\pm$ 1.9053\\
$N=500$&\textbf{4.9135} $\pm$ 0.5328&\textbf{4.9135} $\pm$ 0.5328&4.9145 $\pm$ 0.5328&4.5618 $\pm$ 0.0801&\textbf{0.0455} $\pm$ 0.0033 &15.5149 $\pm$ 1.7017\\
\bottomrule
\end{tabular*}} \vspace{0cm}
\end{table}
\end{minipage}
\end{center}
\cref{table:objsub} compares the average objective values, including the standard deviation after the $\pm$ mark, and computation times of different portfolio optimization methods. Results show that both DP and PAV algorithms outperform \texttt{fmincon} in terms of objective values, with PAV being the fastest method.  The findings confirm the theoretical complexity analysis, indicating that PAV has $O(N)$ time complexity. DP and PAV algorithms always produce the same objective values, indicating that PAV identifies the global minimum like DP. However, non-global minimum points are observed with real data in ADMM subproblems in the next subsection.

\subsection{Numerical Tests for the ADMM  Algorithms}\label{sec:5.2}
In this subsection, we compare the performance of our ADMM algorithm based on the two $\y$-subproblem solvers and the \texttt{fmincon} solver in MATLAB for solving the CPT portfolio optimization
problem  on two datasets. The first is the Fama French 48 Industries (FF48), while the other focuses on 458 risky assets listed in the Standard \& Poor's 500 (S\&P 500) index. Particularly, we adopt the utility function \eqref{eq:sshapeu} with the same parameters in \cref{sec:5.1}, and $\cX=\{\x:~\x\ge0,~\sum_{i=1}^d x_i=1\}$, except that we consider two cases of $B$, i.e., $B=0$ or $B=r_f$, where $r_f$ is the risk-free rate. Note that two datasets have different $r_f$s; the specific values of $r_f$ can be found in Appendix \ref{app:expsetting}.
In the two ADMM algorithms, the $\x$-subproblem was solved using the off-the-shelf solver \textsf{Gurobi} 11.0.0 \citep{gurobi}, and the $\y$-subproblem was solved by either the DP or PAV algorithm in \cref{sec:4}. We refer to the ADMM algorithms using the DP and PAV algorithms as ADMM-DP and ADMM-PAV, respectively.  The details of problem settings are presented in Appendix \ref{app:expsetting}.

\begin{center}
\begin{minipage}{\linewidth}
\begin{table}[H]  \caption{Final objective value and time comparison for CPT optimization problem \eqref{pb:1} with $B=0$. Numbers in parentheses represent iteration numbers of ADMM algorithms. }
\label{table:objcmp}
{\scriptsize
\begin{tabular*}{\textwidth}{l@{\extracolsep{\fill}}|lll|lll}
\toprule
\multirow{2}{*}{Scenarios } &\multicolumn{3}{c|}{objective value }&\multicolumn{2}{c}{time (seconds)}\\
  &ADMM-DP  & ADMM-PAV & \texttt{fmincon}  & ADMM-DP & ADMM-PAV  &\texttt{fmincon} \\
\midrule
$d=$48,$N=$50 &\textbf{-1.855$\times 10^{-3}$} & \textbf{-1.855$\times 10^{-3}$} & -1.364$\times 10^{-3}$
&6.169 (59)& 2.018 (61)& \textbf{1.292}\\
$d=$48,$N=$100 &-3.472$\times 10^{-4}$ & \textbf{-3.477$\times 10^{-4}$} & -1.580$\times 10^{-4}$
&13.37 (62)& \textbf{2.356} (66)& 2.853\\
$d=$48,$N=$150 &\textbf{4.410$\times 10^{-4}$ }& \textbf{4.410$\times 10^{-4}$ }& 1.500$\times 10^{-3}$
&27.19 (63)& 2.679 (68)& \textbf{0.530}\\
$d=$48,$N=$200 &\textbf{6.390$\times 10^{-4}$ }& \textbf{6.390$\times 10^{-4}$ }& 6.537$\times 10^{-4}$
&45.42 (63)& 3.273 (64)& \textbf{1.981}\\
$d=$48,$N=$250 &\textbf{1.196$\times 10^{-3}$ }& 1.199$\times 10^{-3}$ & 1.261$\times 10^{-3}$
&69.4 (65)& 4.973 (61)& \textbf{2.414}\\
$d=$48,$N=$300 &\textbf{2.323$\times 10^{-3}$ }& \textbf{2.323$\times 10^{-3}$ }& 2.400$\times 10^{-3}$
&96.98 (62)& 4.359 (64)& \textbf{0.949}\\
\midrule
$d=$458,$N=$300 &\textbf{3.329$\times 10^{-3}$} & 3.363$\times 10^{-3}$ & 6.828$\times 10^{-3}$
& 112.7 (62)& 20.93 (59)& \textbf{8.743}\\
$d=$458,$N=$400 &\textbf{3.111$\times 10^{-3}$} & \textbf{3.111$\times 10^{-3}$} & 4.093$\times 10^{-3}$
& 192.7 (62)& \textbf{26.52} (61)& 41.17\\
$d=$458,$N=$500 &\textbf{4.605$\times 10^{-3}$} & \textbf{4.605$\times 10^{-3}$} & 4.810$\times 10^{-3}$
& 260.9 (57)& \textbf{24.00} (53)& 583.7\\
$d=$458,$N=$600 &4.995$\times 10^{-3}$ & \textbf{4.994$\times 10^{-3}$} & 5.179$\times 10^{-3}$
& 376.2 (57)& \textbf{24.05} (52)& 325.2\\
$d=$458,$N=$700 &4\textbf{.796$\times 10^{-3}$ }& \textbf{4.796$\times 10^{-3}$} & 4.952$\times 10^{-3}$
& 499.8 (57)& \textbf{26.94} (57)& 378.1\\
$d=$458,$N=$800 &\textbf{4.862$\times 10^{-3}$} & \textbf{4.862$\times 10^{-3}$} & 5.015$\times 10^{-3}$
& 572.3 (48)& \textbf{22.82} (47)& 691.9\\
$d=$458,$N=$900 &\textbf{4.606$\times 10^{-3}$} & \textbf{4.606$\times 10^{-3}$} & 4.759$\times 10^{-3}$
& 771.2 (51)& \textbf{25.10} (51)& 528.5\\
$d=$458,$N=$1000 &\textbf{5.079$\times 10^{-3}$} & 5.083$\times 10^{-3}$ & 5.196$\times 10^{-3}$
& 1051 (56)& \textbf{26.66} (53)& 662.6\\
\bottomrule
\end{tabular*}\vspace{0.2cm}
}
 \vspace{0cm}
\end{table}
\end{minipage}
\end{center}

\begin{center}
\begin{minipage}{\linewidth}
\begin{table}[H]  \caption{Final objective value and time comparison for CPT optimization problem \eqref{pb:1} with $B=r_f$. Numbers in parentheses represent iteration numbers of ADMM algorithms.}
\label{table:objcmp B riskfree}
{\scriptsize
\begin{tabular*}{\textwidth}{l@{\extracolsep{\fill}}|lll|lll}
\toprule
\multirow{2}{*}{Scenarios } &\multicolumn{3}{c|}{objective value }&\multicolumn{2}{c}{time (seconds)}\\
 &ADMM-DP  & ADMM-PAV & \texttt{fmincon}  & ADMM-DP & ADMM-PAV  &\texttt{fmincon} \\
\midrule
$d=$48,$N=$50 &\textbf{-1.777 $\times 10^{-3}$} & \textbf{-1.777 $\times 10^{-3}$} & -1.288 $\times 10^{-3}$
&6.447 (61) &1.963 (61) &0.695\\
$d=$48,$N=$100 &-8.442 $\times 10^{-5}$ & \textbf{-8.873 $\times 10^{-5}$} & -6.573 $\times 10^{-5}$
&15.09 (66) &2.436 (64) &\textbf{0.758}\\
$d=$48,$N=$150 &\textbf{5.238 $\times 10^{-4}$} & \textbf{5.238 $\times 10^{-4}$} & 1.581 $\times 10^{-3}$
&24.90 (58) &2.409 (58) &\textbf{0.453}\\
$d=$48,$N=$200 &7.249 $\times 10^{-4}$ & \textbf{7.248 $\times 10^{-4}$} & 7.419 $\times 10^{-4}$
&36.36 (59) &2.933 (60) &\textbf{1.839}\\
$d=$48,$N=$250 &\textbf{1.282 $\times 10^{-3}$} & \textbf{1.282 $\times 10^{-3}$} & 1.349 $\times 10^{-3}$
&53.50 (61) &3.576 (63) &\textbf{2.403}\\
$d=$48,$N=$300 &\textbf{2.409 $\times 10^{-3}$} & \textbf{2.409 $\times 10^{-3}$} & 2.425 $\times 10^{-3}$
&89.14 (67) &3.379 (58) &\textbf{0.902}\\
\midrule
$d=$458,$N=$300 &\textbf{3.335 $\times 10^{-3}$} & 3.369 $\times 10^{-3}$ & 6.834 $\times 10^{-3}$
&94.98 (57) &19.19 (57)&\textbf{11.65}\\
$d=$458,$N=$400 &3.117 $\times 10^{-3}$ & \textbf{3.115 $\times 10^{-3}$} & 4.109 $\times 10^{-3}$
&161.6 (59) &\textbf{23.06} (58)&32.33\\
$d=$458,$N=$500 &\textbf{4.611 $\times 10^{-3}$} & \textbf{4.611 $\times 10^{-3}$} & 4.814 $\times 10^{-3}$
&205.8 (52) &\textbf{22.79} (53)&632.0\\
$d=$458,$N=$600 &5.001 $\times 10^{-3}$ & \textbf{5.000 $\times 10^{-3}$} & 5.182 $\times 10^{-3}$
&296.5 (52) &\textbf{25.03} (57)&457.7\\
$d=$458,$N=$700 &\textbf{4.802 $\times 10^{-3}$} & 4.803 $\times 10^{-3}$ & 4.961 $\times 10^{-3}$
&407.7 (54) &\textbf{26.21} (56)&701.6\\
$d=$458,$N=$800 &\textbf{4.868 $\times 10^{-3}$} & \textbf{4.868 $\times 10^{-3}$} & 5.019 $\times 10^{-3}$
&488.6 (52) &\textbf{21.48} (47)&124.2\\
$d=$458,$N=$900 &\textbf{4.612 $\times 10^{-3}$} & \textbf{4.612 $\times 10^{-3}$} & 4.764 $\times 10^{-3}$
&544.0 (47) &\textbf{22.03} (47)&224.4\\
$d=$458,$N=$1000 &\textbf{5.086 $\times 10^{-3}$} & 5.092 $\times 10^{-3}$ & 5.201 $\times 10^{-3}$
&770.1 (53) &\textbf{26.68} (56)&577.7\\
\bottomrule
\end{tabular*}\vspace{0.2cm}
}
 \vspace{0cm}
\end{table}
\end{minipage}
\end{center}

In \cref{table:objcmp} (resp. \cref{table:objcmp B riskfree}), we present a comparison of final objective values and CPU times for various algorithms across different numbers of scenarios for the case of $B=0$ (resp. $B=r_f$).
The objective value is the negative of the CPT utility function, with lower values indicating better solutions.
The results on both tables demonstrate that, on average, ADMM-PAV exhibits the best performance. Specifically, in both FF48 instances with $d=48$ and S\&P 500 instances with $d=458$, ADMM-PAV and ADMM-DP achieve much better objective values than \texttt{fmincon}, but do not dominate each other.
Among these, ADMM-PAV achieves the shortest runtime,  at most twenty to thirty seconds in all scenarios, and is the fastest algorithm for cases with $N \geq 400$. In contrast, ADMM-DP is much slower in all instances compared to ADMM-PAV, with its runtime increasingly slowing down as $N$ increases. 

\subsection{Numerical Tests for Comparison with \cite{luxenberg2024portfolio}}\label{sec:5.3}
 In this subsection, we compare our proposed ADMM method with three methods outlined in \cite{luxenberg2024portfolio}.
Firstly, we provide brief descriptions of these methods for context:
\begin{itemize}
    \item \textbf{Minorization-Maximization (MM)}: This method iteratively constructs a concave approximation to the CPT utility function. It then maximizes this approximation to determine the subsequent iterate.
    \item \textbf{Iterated Convex-Concave (CC)}: This method involves fixing the probability weights at each iteration, followed by applying the convex-concave procedure to solve the optimization problem of the CPT utility with fixed weights.
    \item \textbf{Projected Gradient Ascent (GA)}: This method employs an automatic differentiation package to compute the gradient of the objective function at differentiable points and to identify a suitable surrogate for the gradient at nondifferentiable points. Each iteration involves using gradient ascent combined with projection onto the set of portfolio constraints to advance to the next point.
\end{itemize}

We remark that both MM and CC solve a different approximate model, where monotonicity is enforced on the weights $c_i$ to ensure convexity or concavity, rather than \eqref{pb:1}. Note that the CPT utility is not differentiable,  so there is no theoretic guarantee of the projected gradient method, which uses automatic differentiation to compute the ``gradient".
Moreover, one can also show that the convexity properties in Section 2.3 in \cite{luxenberg2024portfolio} fail for the S-shaped power utility function given by \eqref{eq:sshapeu},
thereby hindering the applicability of both MM and CC.
Nevertheless, we conduct a comparison based on the approximate CPT model in \cite{luxenberg2024portfolio}, which imposes the monotonicity of weights and employs an exponential utility function.
For a detailed description of the setup, please refer to Appendix \ref{sec:A5.3}.

\begin{center}
\begin{minipage}{\linewidth}
\begin{table}[H]  \caption{Final objective value and time comparison for CPT optimization problem under the setting of  \cite{luxenberg2024portfolio}. {``\texttimes" indicates that CC fails to solve the current scenarios.}}
\label{table:ccmmgacmp}
{\setlength{\tabcolsep}{1pt}
\scriptsize
\begin{tabular*}{\textwidth}{l@{\extracolsep{\fill}}|cccc|cccc}
\toprule
\multirow{2}{*}{Scenarios } &\multicolumn{4}{c|}{objective value }&\multicolumn{4}{c}{time (seconds)}\\
  &ADMM-PAV  &MM &CC &GA &ADMM-PAV  &MM &CC &GA   \\ \midrule
  $d=$48,$N=$50 &\textbf{-1.954$\times 10^{-2}$} & -1.946$\times 10^{-2}$ & \texttimes & -1.938$\times 10^{-2}$
&\textbf{4.385}  &5.889& \texttimes &182.9\\
$d=$48,$N=$100 &\textbf{-1.030$\times 10^{-2}$} & -1.013$\times 10^{-2}$ & -1.028$\times 10^{-2}$ & -7.943$\times 10^{-3}$
&1.443 &20.32&\textbf{1.058}&185.1\\
$d=$48,$N=$150 &-8.453$\times 10^{-3}$ & -8.367$\times 10^{-3}$ & \textbf{-8.456$\times 10^{-3}$} & -7.267$\times 10^{-3}$
&1.967 &63.25&\textbf{1.398}&190.9\\
$d=$48,$N=$200 &-4.616$\times 10^{-3}$ & \textbf{-6.028$\times 10^{-3}$} & -4.604$\times 10^{-3}$ & -2.548$\times 10^{-3}$
&1.939  &182.4&\textbf{1.477}&195.2\\
$d=$48,$N=$250 &\textbf{-4.877$\times 10^{-3}$} & -4.742$\times 10^{-3}$ & -4.724$\times 10^{-3}$ & -2.704$\times 10^{-3}$
&2.041 &100.1&\textbf{1.781}&200.1\\
$d=$48,$N=$300 &\textbf{-3.726$\times 10^{-3}$} & -3.663$\times 10^{-3}$ & -3.714$\times 10^{-3}$ & -2.584$\times 10^{-3}$
&\textbf{1.525}  &739.9&2.195&204.6\\
\midrule
$d=$458,$N=$300 &\textbf{-1.840$\times 10^{-2}$} & -1.836$\times 10^{-2}$ & \texttimes & 4.269$\times 10^{-3}$
&\textbf{12.16} &1370&\texttimes &209.9\\
$d=$458,$N=$400 &\textbf{-1.678$\times 10^{-2}$} & -1.662$\times 10^{-2}$ & \texttimes & 3.469$\times 10^{-3}$
&\textbf{21.27} &1860&\texttimes &218.6\\
$d=$458,$N=$500 &\textbf{-1.045}$\times 10^{-2}$ & -1.033$\times 10^{-2}$ & -1.039$\times 10^{-2}$ & 8.810$\times 10^{-3}$
&\textbf{21.51}&3693&98.11&228.8\\
$d=$458,$N=$600 &\textbf{-4.343$\times 10^{-3}$} & -3.679$\times 10^{-3}$ & -4.148$\times 10^{-3}$ & 8.575$\times 10^{-3}$
&\textbf{21.09} &2808&101.5&238.5 \\
$d=$458,$N=$700 &\textbf{-5.028$\times 10^{-3}$} & -4.755$\times 10^{-3}$ & -5.008$\times 10^{-3}$ & 8.520$\times 10^{-3}$
&\textbf{17.83}&11590&125.0&247.3\\
$d=$458,$N=$800 &\textbf{-3.384$\times 10^{-3}$} & -3.099$\times 10^{-3}$ & \texttimes & 7.738$\times 10^{-3}$
&\textbf{18.13} &7389&\texttimes &259.5\\
$d=$458,$N=$900 &\textbf{-3.599$\times 10^{-3}$} & -3.340$\times 10^{-3}$ & -3.561$\times 10^{-3}$ & 6.678$\times 10^{-3}$
&\textbf{20.07} &10380&150.6&269.0\\
$d=$458,$N=$1000 &-3.414$\times 10^{-3}$ & -3.014$\times 10^{-3}$ & \textbf{-3.418$\times 10^{-3}$} & 6.613$\times 10^{-3}$
&\textbf{20.58}&18930&179.3&278.5\\
\bottomrule
\end{tabular*}\vspace{0.2cm}
}
 \vspace{0cm}
\end{table}
\end{minipage}
\end{center}
 In \cref{table:ccmmgacmp}, we provide a comparative analysis of the final objective values and CPU times achieved by different algorithms across a variety of scenario counts. Given that ADMM-DP typically underperforms relative to ADMM-PAV, we have opted not to include ADMM-DP in the experimental evaluation within this subsection.

The results also demonstrate that, on average, ADMM-PAV performs the best. The table presented reveals that ADMM-PAV achieves superior objective function values compared to other algorithms in nearly all scenarios, and are significantly faster than the other three methods. Although MM significantly outperforms other algorithms in terms of the objective function value for FF48 instances with \(d=48, N=200\), a closer examination reveals that this is merely an isolated incident. It is highly probable  that the MM algorithm converges to a different local optimal point in this scenario, diverging from the solutions attained by other algorithms. In most cases, its performance is unsatisfactory. Furthermore, in terms of computational time, it is considerably less efficient than the other algorithms. On the other hand, while the CC algorithm does exhibit superior objective function values compared to ADMM-PAV in certain scenarios for both FF48 instances with \(d=48\) and S\&P 500 instances with \(d=458\), we point out that it can fail in some scenarios. This is due to occasional failures in solving subproblems during the iterative process using the ECOS solver, thereby compromising its overall robustness. Regarding the GA algorithm, it can be observed that, apart from achieving results close to those of ADMM-PAV in scenarios with very low dimensions, its performance is poor in other cases.

\subsection{Empirical Study}
We extend our investigation by conducting an empirical study of portfolio optimization under the CPT framework using the proposed method. This study has not been previously conducted due to the absence of solvers. Through this empirical study, we investigate the effects of CPT's parameters on optimal portfolios, providing valuable insights into the roles of CPT's parameters. Detailed results of our empirical study can be found in Appendix \ref{app:es}.

\section{Conclusion}\label{sec:7}

In this paper, we propose a novel numerical method for solving the CPT utility maximization model under constraints with finite realizations of scenarios. Our method employs an ADMM algorithm to handle the challenges posed by the S-shaped utility and inverse-S-shaped probability distortion function of the CPT optimization problem. One of the subproblems in the ADMM is simple, while the other is difficult, which minimizes the negative of the expected utility subject to a chain constraint.
We further propose a DP method for finding the global optimal solution of one subproblem in ADMM. However, due to the high computational cost of the DP method, we also introduce a PAV algorithm that returns a stationary point of a reformulation of the subproblem. Additionally, we conduct comprehensive convergence results for the ADMM algorithm and subproblem solvers. { Our experimental results demonstrate that the proposed ADMM algorithm with the PAV algorithm outperforms both the ADMM algorithm with the DP method, the MATLAB solver \texttt{fmincon} and the methods proposed in \cite{luxenberg2024portfolio}.}



\bibliographystyle{informs2014} 
\bibliography{mybib2.bib} 

%
%
%


\appendix
\section{Proofs}\label{ap:proof}
\subsection{Proof of \cref{lem:locallip}}\label{proof:l1}
\begin{myproof}
As $\frac{\sigma}{2} \left \|\y - R\x + \frac{\blambda}{\sigma} \right\|^2$ is differentiable, we only need to show that $\Omega(\y)$ is locally Lipschitz.
Let  $\mathcal B(\bar\y,\epsilon)$ denote the open ball with radius $\epsilon$ centered at $\bar\y$. 
First note that $U(\theta)$ is locally Lipschitz for $\theta\neq B$ as $U$ is differentiable at $\theta\neq B$.
When $\theta = B$, if $ U'(B) < \infty$, $U(\theta)$ is also locally Lipschitz at $B$ as its left-derivative and right-derivative are bounded. If $y_1=\dots=y_N$, then it is easy to verify that $\Omega(\y)$ is locally Lipschitz as $U'(y_i)<\infty,~\forall i$ for both cases of $y_i\neq B$ and  $y_i=B$.

Now let us consider the case that not all $y_i$ are equivalent. Let
$\epsilon=\min_{i\in I}\{\bar y_{[i+1]}-\bar y_{[i]}\},$
where $I=\{i:y_{[i+1]}>y_{[i]},~i=1,\ldots,N-1\}$.
It is obvious that $\epsilon>0$.
Let $\y^1,\y^2\in\mathcal B(\bar\y,\epsilon)$ be  two arbitrary points. Thanks to the definition of $\epsilon$, there exist two permutations $\{j_1,j_2,\ldots,j_N\}$ and $\{k_1,k_2,\ldots,k_N\}$ of $\{1,2,\ldots,N\}$ such that
 \[y^1_{j_i}\le y^1_{j_{i+1}}, \quad \bar y_{j_i}\le \bar y_{j_{i+1}}\]
  and
 \[y^2_{k_i}\le y^2_{k_{i+1}}, \quad \bar y_{k_i}\le \bar y_{k_{i+1}}\]
 for $i=1,2,\ldots,N$. (Note that the two permutations may be different if $\bar y_{[i]}=\bar y_{[i+1]}$ for some $i$.)
Therefore we have
\begin{equation*}
\begin{aligned}
 |\Omega(\bar\y)-\Omega(\y^1)|
 & =  |\sum_{i = 1}^{N} c_i U(\bar y_{j_i}) - c_i U(y^1_{j_i})|\\
  & \le \sum_{i = 1}^{N} c_i |U(\bar y_{j_i}) - U(y^1_{j_i})|\\
 &\le \max_{1\le i \le N}{c_i} L_{j_i}|\bar y_{j_i}-y^1_{j_i}| \\
 &\le L^1\|\bar\y-\y^1\|,
\end{aligned}
\end{equation*}
where $L_{j_i}$ is the Lipschitz constant for $U$ in $B(\bar y_{j_i},\epsilon)$ and $L^1=\max_{1\le i \le N}{c_i} L_{j_i}$.
Similarly, we also have
\[
 |\Omega(\bar\y)-\Omega(\y^2)|\le L^2\|\bar\y-\y^2\|
\]
for some constant $L^2$.
This implies that
\[
 |\Omega(\y^1)-\Omega(\y^2)|\le  |\Omega(\y^2)-\Omega(\bar\y)| + |\Omega(\y^1)-\Omega(\bar\y)|\le (L^1+L^2)\|\y^1-\y^2\|.
\]
This completes the proof.
\end{myproof}

Our proof is fundamental and simple. The proof can also be used to show the Locally Lipschitz property for the empirical VaR in \cite{cui2018portfolio}, which can be represented as $\sum_{i=1}^Na_i\y_{[i]}$ for some $a_i\ge0$, admitting a simpler form than  $\Omega(\y)$.

\subsection{Proof of \cref{thm:admm}}\label{ap:admm}
\begin{myproof}
Note that in the $(k+1)$-th iteration, $\x^{k+1}$ minimizes $L_{\sigma}(\x,\y^k;\blambda^k) $.    Because $L_{\sigma}(\x,\y^k;\blambda^k) $ is a strongly convex function of $\x$, we have
\begin{equation}\label{xx1}
        L_{\sigma}(\x^k,\y^k;\blambda^k) - L_{\sigma}(\x^{k+1},\y^k;\blambda^k) \geq \nabla_{\x} L_{\sigma}(\x^{k+1},\y^k;\blambda^k)^T (\x^k - \x^{k+1}) + \frac{\sigma\rho}{2}\| \x^k - \x^{k+1} \|^2,
\end{equation}
where $\rho>0$ is the minimal eigenvalue of $R^TR$.  Because $\x^{k+1}$ is the global minimum of the $\x$-subproblem and $\x_k\in\mathcal X$,  we must have
\begin{equation}\label{xx2}
        \nabla_{\x} L_{\sigma}(\x^{k+1},\y^k;\blambda^k)^T (\x^k - \x^{k+1}) \geq 0.
\end{equation}
Combining \eqref{xx1} and \eqref{xx2} gives
\begin{equation}\label{xxx}
        L_{\sigma}(\x^k,\y^k;\blambda^k) - L_{\sigma}(\x^{k+1},\y^k;\blambda^k) \geq \frac{\sigma\rho}{2}\| \x^k - \x^{k+1} \|^2.
\end{equation}

Note that the update \eqref{eq:sublam} is equivalent to
\begin{equation}\label{boundy}
        \y^{k+1} - R\x^{k+1} = - \frac{1}{\sigma}(\blambda^{k} - \blambda^{k+1}).
\end{equation}
Thus we have
\begin{equation}\label{lll}
        L_{\sigma}(\x^{k+1},\y^{k+1};\blambda^k) - L_{\sigma}(\x^{k+1},\y^{k+1};\blambda^{k+1}) = \langle \blambda^{k} - \blambda^{k+1} , \y^{k+1} - R\x^{k+1} \rangle = -\frac{1}{\sigma}\| \blambda^{k} - \blambda^{k+1} \|^2.
\end{equation}
Summing up \eqref{xxx}, \eqref{eq:ydescent} and \eqref{lll} gives rise to
\begin{equation}\label{iter}
        L_{\sigma}(\x^k,\y^k;\blambda^k) - L_{\sigma}(\x^{k+1},\y^{k+1};\blambda^{k+1}) \geq
\frac{\sigma\rho}{2} \| \x^k - \x^{k+1} \|^2 -\frac{1}{\sigma}\| \blambda^{k} - \blambda^{k+1} \|^2.
\end{equation}
Under the assumption that the sequence $\{ \blambda^k \}$ is bounded and the set $\cX$ is compact (implying  $\{ \x^k \}$ is bounded), it follows from \eqref{boundy} that $\{\y^k\}$ is also bounded.
Therefore $L_{\sigma}(\x^k,\y^k;\blambda^k)$  is also bounded.
Then summing  up both sides of \eqref{iter} for $k$ from 1 to infinity, we have
\begin{equation*}
        \frac{\sigma\rho}{2} \sum_{k=1}^{\infty}\| \x^{k+1} - \x^k\|^2 < \infty + \frac{1}{\sigma}\sum_{k = 1}^{\infty} \| \blambda^{k+1} - \blambda^k \|^2.
\end{equation*}
This, together with the assumption that $\sum_{k = 1}^{\infty} \| \blambda^{k+1} - \blambda^k \|^2 < \infty$, implies $\sum_{k=1}^{\infty}\| \x^{k+1} - \x^k\|^2 < \infty$, which further yields $\| \x^{k+1} - \x^k\| \rightarrow 0$ for $k \rightarrow \infty $.
Using $\sum_{k = 1}^{\infty} \| \blambda^{k+1} - \blambda^k \|^2 < \infty$ again,  we have $\blambda^{k+1}-\blambda^k\to\bz$.
Then \eqref{boundy} further implies
\begin{equation}\label{eq:yRx}
\y^{k} - R\x^{k} \rightarrow \bz.
\end{equation}
Due to $\| \x^{k+1} - \x^k\| \rightarrow 0$ and $\blambda^{k+1}-\blambda^k\to\bz$,  \eqref{boundy} implies that
\begin{equation}\label{eq:ydiffto0}
\|\y^{k+1}-\y^k\|\le \|R(\x^{k+1}-\x^k)\|+\frac{1}{\sigma}\|(\blambda^k-{\blambda^{k+1}})-(\blambda^{k-1}-\blambda^k)\|\to 0.
\end{equation}

As the $\x$-subproblem is globally solved, we must have
        \begin{align}
                 \bz &\in - \sigma R^T\left(\y^k - R\x^{k+1} + \frac{\blambda^k}{\sigma}\right) +  N_{\mathcal X}(\x^{k+1}). \label{xkfirst}
        \end{align}
Note also that \eqref{ykfirst} implies
\begin{equation}\label{yre2}
        \bz \in R^T \partial\Omega(\y^{k+1}) + \sigma R^T\left(\y^k - R\x^{k+1} + \frac{\blambda^k}{\sigma}\right) + \sigma R^T(\y^{k+1} - \y^{k}).
\end{equation}
From \eqref{xkfirst} and \eqref{yre2}, we deduce that
\begin{equation*}\label{yre3}
        \bz\in R^T \partial\Omega(\y^{k+1}) + N_{\mathcal X}(\x^{k+1}) + \sigma R^T\left(\y^{k+1} - \y^{k}\right).
\end{equation*}
This, together with \eqref{eq:ydiffto0}, \cref{lem:hemic} and the fact that the normal cone of a closed set is outer semicontinuous  set-valued mapping \citep[Proposition 6.6]{rockafellar2009variational}, yields
\[\bz\in R^T \partial\Omega( \y^*) + N_{\mathcal X}(\x^{*}),\]
where $(\x^*,\y^*)$ is any accumulation point of $\{(\x^k,\y^k)\}$. Note also that
\eqref{eq:yRx} implies
$\y^*=R\x^*$. The proof is completed.
\end{myproof}

 \subsection{Proof of \cref{lem:3ora}}\label{ap:lem3}
\begin{myproof}
Consider the interval $(B,\infty)$. From \cref{fact}, we obtain that it needs at most one root-finding oracle to find the local minimum in $(B,\infty)$ in case that $\zeta'(B^+)<0$.

Consider the interval $(-\infty,B)$.
If $\zeta{''}(B^-) \ge0$,   we obtain that it needs at most one root-finding oracle to find the local minimum in $(w,B)$ in case of $\zeta'(B^-)>0$.
If $\zeta{''}(B^-) < 0$,  according to \cref{fact}, we first use one oracle to find the point $C$ such that $\zeta''(C)=0$, thanks to the monotonicity of $\zeta''(C)$ in $(-\infty,B)$. Then we need further to consider the case $\zeta'(C)>0$. We need one oracle to find the local minimum in $(w,C)$.

To identify if $B$ is a local minimum, we only need to check $\zeta'(B^-)$ and $\zeta'(B^+)$ without a root-finding oracle.
\end{myproof}

\subsection{Proof of \cref{lem:4}}\label{ap:lem4}
\begin{myproof}
Suppose for some index $i$ we have $U'(y^*_i)= \infty$, in which case we must have $y^*_i=B$.
It is straightforward to verify that
\begin{equation*}
\lim_{y\to B^{-}}f'_i(y)= \lim_{y \to B^+}f'_i(y)  =   -\infty,\quad \forall i.
\end{equation*}
Then we can perturb $\y^*$ a bit by slightly enlarging all $y^*_j$ such that $y^*_j=y^*_i$ and do not violate the chain constraint.
This yields a smaller objective value and still satisfies the chain constraint, which contradicts the optimality of $\y^*$.
\end{myproof}

 \subsection{Proof of \cref{lem:bound}}\label{ap:bound}
\begin{myproof}
First note that $f_i(\cdot)$ is in the form of $\zeta(\cdot)$.
Note also that $ y_1^*\le \dots \le y_N^* $ due to the chain constraint in  \eqref{eq:dpyf}.
Suppose on the contrary that $y_1^*\le\cdots\le y_k^* < w_1\le y_{k+1}^*$.
Define $\tilde{\y}$ as $\tilde{y}_i = \max(w_1, y^{*}_i)$, $i=1,\ldots,N$, which is also feasible to \eqref{eq:dpyf}.
It follows from \cref{asmp:1} and \eqref{d1} that $f_i{'}(u) < 0$ for all $u\in(-\infty,w_1)$, $i=1,\ldots,k$.
This yields $\sum_{i=1}^N f_i(y_i^*) > \sum_{i=1}^N f_i(\tilde{y}_i)$, which contradicts the optimality of $\y^*$. Hence we must have $
w_1\le y_1^*\le \dots \le y_N^*$.

 Thanks to the fact that $f_i'(u)>0$ for all $u\in( \max \left( B+1, w_i + \frac{b_iU'(B+1)}{\sigma}\right),\infty)$, due to our discussion of $\zeta(y)$ in \cref{fact}, the upper bound can be analyzed similarly
\end{myproof}

\subsection{Proof for \cref{lem:recur}}\label{proof:l7}
\begin{myproof}
We prove the lemma using induction.
First note that $h_1$ can be expressed in \eqref{eq:hexpre} as $h_1(z) = 0$ by $a_k^{m},b_k^{m},\sigma_k^{m} = 0$. Next we assume for $m=n-1$, \eqref{eq:hexpre} holds.

Now let us show \eqref{eq:hexpre} holds for $m=n$.
It is straightforward from \eqref{eq:phiy} and the fact that $f_i$ is of expression \eqref{eq:zeta} that $\phi_{n-1}(z)$ has the same expression form as \eqref{eq:hexpre},
\begin{equation*}
        \phi_{n-1}(z) = \zeta_{\hat{a}_{n-1,k},\hat{b}_{n-1,k},\hat{w}_{n-1,k},\hat{\sigma}_{n-1,k},U}(z)+ \hat{M}_{n-1,k},\ z\in  [s^{n-1}_k,s^{n-1}_{k+1}],
\end{equation*}
where $\hat{a}_{n-1,k} = a_{n-1,k} + a_{n-1}$, $\hat{b}_{n-1,k} = b_{n-1,k} + b_{n-1}$,
        $\hat{\sigma}_{n-1,k}  = \sigma_{n-1,k} + \sigma$,
        $\hat{w}_{n-1,k}  = \frac{\sigma_{n-1,k}w_{n-1,k}+ \sigma w_{n-1}}{\sigma_{n-1,k} + \sigma}$, and
        $\hat{M}_{n-1,k}  = M_{n-1,k} +\frac{\sigma_{n-1,k}w_{n-1,k}^2 + \sigma w_{n-1}}{2} - \frac{(\sigma_{n-1,k}w_{n-1,k}+ \sigma w_{n-1})^2}{2(\sigma_{n-1,k} + \sigma)}$.
Hence $\phi_{n-1}(z)$ is in the form of  $\zeta (y)$ plus a constant.
Let \[T^{n-1}_{k} = \Call{Decompose}{[s^{n-1}_k,s^{n-1}_{k+1}], \phi_{n-1}(\cdot)}.\] For every smaller interval $[t_i,t_{i+1}]$ in $T^{n-1}_{k}$, we obtain that $\phi_{n-1}$ is either strictly increasing or strictly decreasing due to \cref{def:decom}.

Then we deduce $h_n(z)$ sequentially for every $[t_i,t_{i+1}]$ in $T^{n-1}_{k}$ using the \Call{Update}{}  procedure.  Note that $h_n(t_i)$ is known when we proceed in  $[t_i,t_{i+1}]$ as a result of the previous interval.
When $t_i = s^{n-1}_1 = l_b$, which is the very beginning of all intervals, we assign $h_n(t_i) = \phi_{n-1}(l_b)$ because $\phi_{n-1}$ is strictly decreasing in $(-\infty,l_b)$. The latter is because $\phi_{n-1}$ is formed by the sum of several $f_i$s and  each $f_i$ is strictly decreasing in $(-\infty,l_b)$ according to the proof of \cref{lem:bound}.
It follows from \eqref{eq:hleqphi} that $h_n(t_i) \leq \phi_{n-1}(t_i)$. Then using $\Call{Update}{[t_i,t_{i+1}], \phi_{n-1}(\cdot), h_n(t_i)}$, we can obtain the expression of $h_n$ in $[t_i,t_{i+1}]$.
Specifically, from the \Call{Update}{} step we know that $h_n(\cdot)$ is a piecewise function that adopts the form of either $\phi_{n-1}(\cdot)$ or a constant in each piece. Thus $h_n(z)$ takes the form of \eqref{eq:hexpre} for $z\in[t_i,t_{i+1}]$.
By sequentially updating $h_n(z)$ in each interval of $[s^{n-1}_k,s^{n-1}_{k+1}]$, we obtain the explicit formula of $h_n(z)$ for $z \in [l_b,u_b]$.
\end{myproof}

\subsection{Proof for \cref{lem:hnconstant}}\label{proof:l9}
\begin{myproof}
Because we merged consecutive intervals that $h_n$ are identical constants and $h_n$ is a non-increasing function,  the fact that $h_n(y_n)$ is a constant in $[s^{n}_{k},s^{n}_{k+1}]$ implies that in the previous interval $[s^{n}_{k-1},s^{n}_{k}]$, $h_n(y_n)$ is strictly decreasing. From Lemma \ref{hndecrease}, we have $h_n(s^{n}_{k}) = \phi_{n-1}(s^{n}_{k})$. The continuity of $h_n(\cdot)$ further guarantees that $h_n(z) = h_n(s^{n}_{k})=\phi_{n-1}(s^{n}_{k})$ for all  $z\in[s^{n}_{k},s^{n}_{k+1}]$.
From \eqref{hphi}, we immediately have that
$\min_{s_k^{n}\le y\le s_{k+1}^{n}}\phi_{n-1}(y) \ge h_n(s_k^n) = \phi_{n-1}(s^{n}_{k})$, and thus $\phi_{n-1}(y) \geq \phi_{n-1}(s^{n}_{k})$  for all $y\in [s^{n}_{k},s^{n}_{k+1}]$.
\end{myproof}

\subsection{Proof for \cref{thm:3}}\label{proof:thm3}
\begin{myproof}
Let the final blocks of PAV be as $\{[s_1+1,s_2],[s_2+1,s_3],\cdots,[s_k+1,s_{k+1}]\}$, where $s_1=0$ and $s_{k+1}=N$.
According to \cref{alg:pav}, we obtain that
$V_{[s_t+1,s_{t+1}]}$ is a global minimizer of
\[
\min_y \sum_{i=s_t+1}^{s_{t+1}} f_{i}(y),
\]
and thus we have
\begin{equation}\label{eq:deriv}
 0\in \sum_{i=s_t+1}^{s_{t+1}} \partial f_{i}(V_{[s_t+1,s_{t+1}]})
\end{equation}
for $t=1,2,\dots,k$.
Now let $\bar y_{s_t+1}=\dots=\bar y_{s_{t+1}}=V_{[s_t+1,s_{t+1}]}$ for $t=1,\dots,k$, and
\[
g_{i}\in \partial f_i(V_{[s_t+1,s_{t+1}]}) = \partial f_i(\bar y_i) = -c_i\partial U(\bar y_{i})+\sigma(\bar y_i-w_i),~i=s_t+1,\ldots,s_{t+1}\]
be such that  $\sum_{i=s_t+1}^{s_{t+1}}g_i=0$ for $t=1,\dots,k$.
Let $\prod\limits_{i=1}^n A_i$ be the Cartesian product of the sets $A_1, A_2, \cdots, A_n$ which represents that
\[
\prod_{i=1}^n A_i = \{(a_1, a_2, \ldots, a_n) \mid a_1 \in A_1, a_2 \in A_2, \ldots, a_n \in A_n\}.
\]
Then \eqref{eq:deriv} implies that
\[
v_j\in \prod_{i=s_t+1}^{s_{t+1}} - c_{i}\partial U(\bar y_{i})+\sigma(\bar y_{i}-w_i), \quad v_j=(g_j,g_{j+1},\dots,g_{j+s_{t+1}-s_t-1})^T,
\]
for  $j=s_t+1,\dots, s_{t+1}$,
due to $w_1\le w_2\le\dots\le w_N$, \eqref{eq:clarke}, and the definition of $g_i$.
Here we use the convention $g_{j+s_{t+1}-s_t-1}= g_{j-1}$ if $j>s_t+1$, which means that $v_j$ varies cyclically as $j$ increases.
This further implies $\bz\in \partial \Gamma(\bar \y)$.
\end{myproof}

\subsection{Proof for \cref{thm:4}}\label{proof:thm4}
\begin{myproof}
First note that line 2 invokes the oracle at most $3N$ times due to \cref{lem:3ora}.
Next consider the pooling adjacent violators period (lines 4 in \cref{alg:pav}). Once we pool two blocks, the number of blocks in $J$ is reduced by 1. Therefore we need at most $N-1$ merging steps. So line 4 invokes the oracle $3(N-1)$ times, thanks to \cref{lem:3ora}. In summary, \cref{alg:pav} invokes the oracle at most  $6N-3$ times.
\end{myproof}

\section{Settings of Experiments in \Cref{sec:5.2}}\label{app:expsetting}

To compile the FF48 dataset, we retrieve daily returns data from Kenneth R. French’s website (\url{http://mba.tuck.dartmouth.edu/pages/faculty/ken.french/data_library.html}) spanning from December 14th, 2016 to February 23th, 2018. Additionally, we gather daily one-month Treasury bill rates from the U.S. Department of the Treasury's website (\url{https://home.treasury.gov/}) for the same period. The mean of FF48 Industries' daily returns during this period ranges from -0.16\% to 0.12\%, while the standard deviations of daily returns range from 0.68\% to 2.00\%. The average one-month Treasury bill rate, which is 0.0034\% per day, acts as the constant risk-free rate. Throughout the experiment, we utilize historical returns data for the first $N=50,100,150,200,250,300$ days as scenarios for the optimization problem.
	
To obtain the S\&P 500 dataset, we collect historical returns data for the assets listed consistently from August 31th, 2009 to August 20th, 2013. The daily return data is sourced from Wharton Research Data Services (\url{https://wrds-www.wharton.upenn.edu/}). Concurrently, we gather daily one-month Treasury bill rates from the U.S. Department of the Treasury's website during the same period. Similar to the FF48 instance, we use the average one-month Treasury bill rate of 0.00027\% per day as the constant risk-free rate. For this optimization problem, we also select historical returns data for the first $N=300,400,500,600,700,800,900,1000$ days as scenarios.

We follow typical  stopping criteria in ADMM \citep{boyd2011distributed}, where  $\| \y^{k} - R\x^{k} \| < \epsilon_1$  and $\| \y^{k} - \y^{k-1}\| < \epsilon_2$ are referred to as primal and dual feasibility, respectively.

To better facilitate convergence,   we dynamically update the penalty parameter $\sigma$ after updating $\lambda$ as
$
\sigma^{k+1}=\left\{\begin{aligned}\iota \sigma^{k}, &\quad\text{ if } \operatorname{mod}(k,5)=0,\\ \sigma^k, &\quad \text{ otherwise},\end{aligned}\right.
$
where $\iota >1$ is a constant.
In both ADMM algorithms, we set $\sigma^0 = 0.7, ~\iota = 1.7,\ \epsilon_1 = 5\times10^{-5}$, and $\epsilon_2 = 2\times10^{-5}$ for FF48 and  $ \epsilon_2 = 5 \times10^{-5}$ for  S\&P 500. We set maximal iteration number as 1000, maximal time  as one hour, and initialize $\x^0=\e/d$, $\y^0=\blambda^0 = \bz$, where $\e$ is the all-one vector.

In both ADMM algorithms, the objective values are computed directly from $\x$, without utilizing the intermediate variable $\y$.
We also employ the \texttt{fmincon} solver to directly solve problem \eqref{pb:1}, with parameters set identically to those in Section \cref{sec:5.1}.
Note that problem \eqref{pb:1} is nonconvex and nonsmooth, so there is no theoretical guarantee of the solution returned by \texttt{fmincon}.

\section{Settings of Experiments in  \cref{sec:5.3}}\label{sec:A5.3}
Following \cite{luxenberg2024portfolio}, the utility function used in  \cref{sec:5.3} is defined as
\begin{equation}\nonumber
U(z)=\begin{cases}
   -1+\exp(\delta_-(z-B)), & \text{if } z \leq B,  \\
   1-\exp(\delta_+(B-z)), & \text{if } z > B,
\end{cases}
\end{equation}
which is S-shaped. Here, $\delta_-$ and $\delta_+$ are hyperparameters.
The probability weighting function is defined as
$$
w_+(p) = \frac{p^{\gamma_+}}{(p^{\gamma_+}+(1-p)^{\gamma_+})^{1/\gamma_+}}, \qquad
w_-(p) = \frac{p^{\gamma_-}}{(p^{\gamma_-}+(1-p)^{\gamma_-})^{1/\gamma_-}},
$$
where $\gamma_-$ and $\gamma_+$ are the probability distortion parameters for losses and gains, respectively.
The the positive and negative decision weights are given by
\begin{align*}
\pi_{+,j}^{\prime} &= \begin{cases}
w_+((N-j+1)/N)-w_+((N-j)/N), & \text{for } j=1,\ldots,N-1, \\
w_+(1/N), & \text{for } j=N,
\end{cases} \\
\pi_{-,j}^{\prime} &= \begin{cases}
w_-((N-j+1)/N)-w_-((N-j)/N), & \text{for } j=1,\ldots,N-1, \\
w_-(1/N), & \text{for } j=N.
\end{cases}
\end{align*}
Let $m_{+} = \arg\min( \pi_+^{\prime})$ and  $m_{-} = \arg\min( \pi_-^{\prime})$. The adjusted decision weights are defined as
$$
\pi_{+,j} = \begin{cases}
\pi_{+,m_{+}}^{\prime}, & \text{for } j=1,\ldots,m_{+}, \\
\pi_{+,m_{+}}^{\prime}, & \text{for } j=m_++1,\ldots,N,
\end{cases} \\ \qquad
\pi_{-,j} = \begin{cases}
\pi_{-,m_{-}}^{\prime}, & \text{for } j=1,\ldots,m_{-}, \\
\pi_{-,m_{-}}^{\prime}, & \text{for } j=m_-+1,\ldots,N.
\end{cases}
$$
Accordingly, we have $a_i = \pi_{-,N-i+1}$ and $b_i = \pi_{+,i}$, for $i=1,...,N$, as introduced after \eqref{eq:subyiso} in Section \ref{sec:4}.
Then the CPT optimization problem is equivalent to
\[
\min_{\x} \quad -\sum_{i = 1}^N  c_i((R\x)_{[i]})  U((R\x)_{[i]})
  \quad {\rm s.t.}~\x \in \mathcal X,
\]
where $c_i(\cdot)= a_i\mathbbm1_{\{\cdot\le B\}}+ b_i\mathbbm1_{\{\cdot> B\}}$.
The parameters are set as $\delta_-=11.4, \delta_+=8.4, \gamma_-=0.79$, and $\gamma_+=0.77$, consistent with \cite{luxenberg2024portfolio}. The maximum number of iterations is set to 1,000 for the MM and CC methods, and 100,000 for the GA method.
The code for the MM, CC and GA methods are from \cite{luxenberg2024portfolio}, all of which is available in their code repository at \href{https://github.com/cvxgrp/cptopt}{https://github.com/cvxgrp/cptopt}.
The ADMM settings are analogous to those described in the preceding subsection.
All algorithms initiate with the portfolio weights set to their average values.

\section{Empirical Investment Performance}\label{app:es}
 In this section, we evaluate the effectiveness of our proposed portfolio optimization method by conducting an empirical analysis. Specifically, we investigate the impact of CPT's four key features on the optimal portfolio when investing in Fama French 48 Industry Indices of the US market (FF48). To this end, we use historical daily data of the 48 industry indices, covering the period from December 14th, 2016 to December 1st, 2021. Notably, the average daily one-month Treasury bill rate during this period, regarded as the risk-free rate, is 0.0041\% per day, which is very close to zero and distinct from the average return of the 48 industry indices, ranging from -0.016\% to 0.129\%.

Given the characteristics of the data, we just choose the case where $B=0$ as our benchmark model. Following the approach of Tversky and Kahneman (1992), we set the other parameters of our benchmark model as $\mu =2.25$, $\alpha =0.88$, $\delta=0.69$, $\gamma=0.61$. This configuration is referred to as ``benchmark" in Table \ref{table:fivemodel}.
To assess the impact of CPT's four key features on the optimal portfolio, we conduct six experiments by varying one feature at a time while keeping the others fixed. The benchmark serves as the first experiment.
The second experiment,
``risk-free reference point'', uses the risk-free rate 0.0041\% as the reference point, and the third experiment,
``large reference point", uses the average return of 48 industry indices in the dataset (0.072\%) as the reference point, showing the impact of different reference points other than zero on the optimal portfolio.  The fourth experiment, ``no risk aversion and risk-seeking", sets $\alpha = 1$, removing risk aversion in the gain domain and risk-seeking in the loss domain, revealing their impact on the optimal portfolio. The fifth experiment, ``no loss aversion", sets $\mu = 1$, investigating the effect of loss aversion on the optimal portfolio. The sixth experiment, ``no probability distortion", sets $\delta = 1$ and $\gamma = 1$, removing probability weighting and examining its role. By comparing the benchmark with each experiment, we gain insights into how each key feature of CPT influences the composition of the optimal portfolio.

\begin{center}
\begin{minipage}{\linewidth}
\begin{table}[H]\scriptsize  \caption{Six model settings.}\label{table:fivemodel}
	\setlength\tabcolsep{15pt}
	\begin{tabular*}{\textwidth}{c |@{\extracolsep{\fill}} ccccc}
		\toprule
		Setting name & $\mu$ & $\alpha$ &  $\delta$ &  $\gamma$ &  $B$   \\
		\midrule
		benchmark & 2.25 & 0.88 & 0.69 & 0.61 & 0 \\
           { risk-free reference point} &
           2.25 & 0.88 & 0.69 & 0.61 & 0.0041\% \\
		{ large reference point} & 2.25 & 0.88 & 0.69 & 0.61 & 0.072\% \\
		{ no risk aversion and risk seeking} & 2.25 & 1 &0.69& 0.61 & 0 \\
		no loss aversion & 1 & 0.88 &0.69& 0.61 & 0 \\
		no probability distortion & 2.25 & 0.88 &1 & 1 & 0 \\
		\bottomrule
	\end{tabular*}
\end{table}
\end{minipage}
\end{center}


The analysis is conducted using a rolling window scheme. Each day in the out-of-sample period from December 12th, 2017, to December 1st, 2021, marks the end of a window. Within each window, we utilize the most recent 250 daily returns as input scenarios for our CPT portfolio optimization model. By applying the ADMM-PAV method, we compute optimal portfolios for the six models and apply these portfolios to the subsequent day. This process generates 1000 different realizations for each portfolio throughout the out-of-sample period.


To compare the optimal portfolios from different experiments, we evaluate several metrics: the cumulative returns achieved by each portfolio during the out-of-sample period; the performance of each portfolio during the out-of-sample period, considering factors such as the mean of return, the volatility of return, Sharpe ratio and maximum drawdown (MD); the diversification levels of portfolio weights. In our comparison, we also include the performance of a simple equally weighted portfolio for reference. This allows us to gauge the relative effectiveness and efficiency of the CPT-based optimization models compared to a basic allocation strategy.

\begin{figure}[H]
    \centering
    \includegraphics[width=0.9\textwidth]{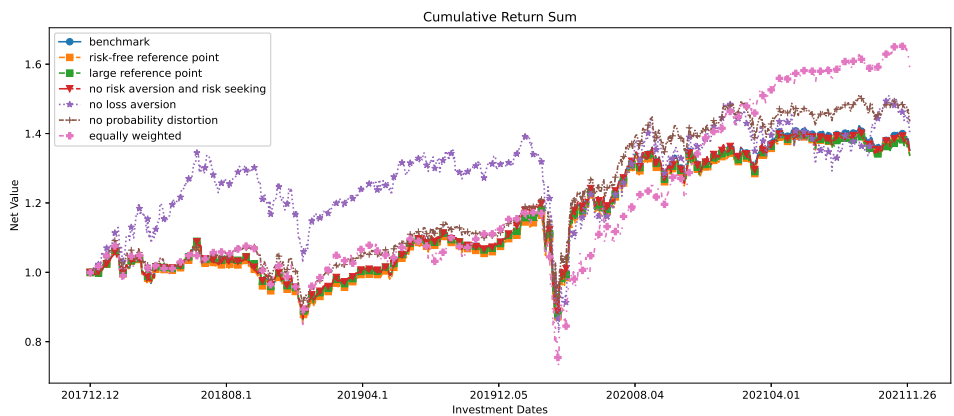}
    \caption{The cumulative returns of seven portfolios.}
    \label{fig:daily}
\end{figure}

We first evaluate the performance of the six optimal portfolios  as well as the equally weighted portfolio based on their cumulative returns over the entire out-of-sample period. Figure \ref{fig:daily} illustrates the cumulative returns of the seven portfolios, with the equally weighted portfolio achieving the highest cumulative return. This is primarily due to the fact that the equal-weighted portfolio is well-diversified but the others are not, as demonstrated by the diversification level analysis provided at the end of the appendix.

Additionally, we calculate several performance metrics of the {\color{blue} seven} portfolios based on the 1000 different realizations. Table \ref{tab:performance} summaries these metrics.

\begin{center}
\begin{minipage}{\linewidth}
\begin{table}[H]\scriptsize
\caption{Performance of seven portfolios.}
 \begin{tabular*}{\textwidth}{@{\extracolsep{\fill}}c|cccc}
\toprule
 Setting & Mean & Volatility & Sharpe Ratio & MD  \\

\midrule
 benchmark& 0.0887 &0.1816&0.4884 &0.2937   \\
 risk-free reference point& 0.0849 &0.1814&0.4682 &0.2934   \\
 large reference point& 0.0845 &0.1823&0.4637 &0.3031   \\
 no risk aversion and risk seeking& 0.0873 &0.1836&0.4753 &0.31   \\
 no loss aversion& 0.1011 &0.2712&0.3729 &0.574   \\
 no probability distortion& 0.1096 &0.1893&0.5789 &0.285   \\
 equally weighted& 0.1483 &0.2213&0.6699 &0.451   \\
  \bottomrule
\end{tabular*}
\label{tab:performance}
\end{table}
\end{minipage}
\end{center}

We conduct an in-depth analysis of the performance of seven portfolios. Transitioning the reference point from 0 to the risk-free rate results in the portfolio displaying a lower mean and Sharpe ratio, while its volatility and MD  remain relatively unchanged. This shift unfavorably impacts assets with small mean and volatility, reducing the feasible set and contributing to the overall decrease in the portfolio's mean and Sharpe ratio.
Similarly, adjusting the reference point from 0 to the average return of $0.072\%$ leads to the portfolio showcasing a lower mean and Sharpe ratio, alongside higher volatility and MD. Investors typically respond by favoring assets with a higher probability of gains to surpass this elevated reference point, despite the associated increase in risk and decrease in overall returns. When we eliminate risk aversion and risk-seeking behavior, the portfolio shows a reduced mean return, heightened volatility and MD, resulting in the lower Sharpe ratio. The removal of loss aversion, however, leads to higher mean returns albeit with increased risk, reflected in the highest volatility and MD and the lowest Sharpe ratio. This aligns with expectations, as loss aversion typically makes investors more risk-averse due to a significant shift in attitude around the reference point. Further analysis without loss aversion indicates that an investor without this bias tends to be less risk-averse and may opt for the riskiest portfolio among the available options. Eliminating probability distortion leads to the portfolio achieving the second-highest Sharpe ratio, possibly due to improved diversification resulting from this adjustment. Our diversification analysis supports these conclusions, showing that the equally weighted portfolio, which is well-diversified, achieves the highest mean return, the second-highest risk, and the highest Sharpe ratio among the portfolios.

\begin{figure}[h]
        \centering

\subfigure[benchmark]{
        \begin{minipage}[t]{0.3\textwidth}
                \centering
                \includegraphics[width=2in]{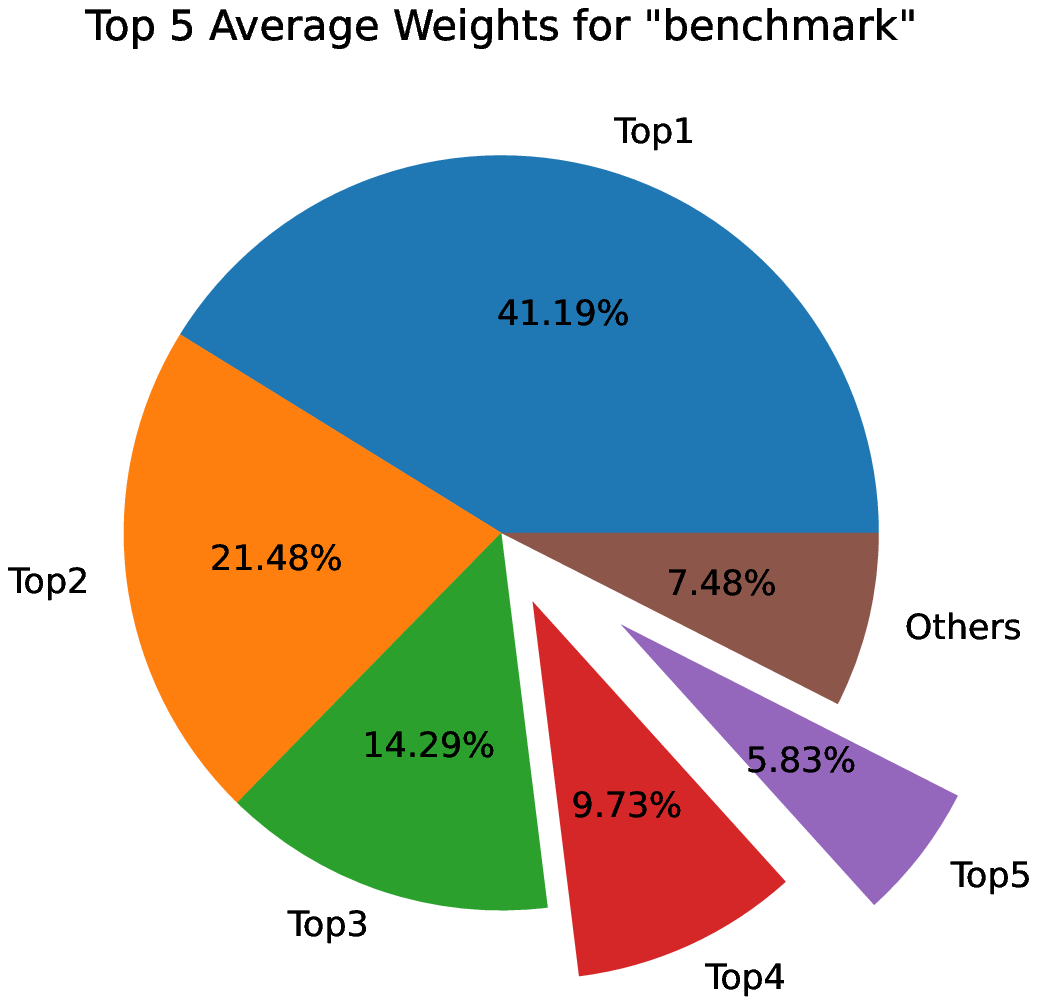}
        \end{minipage}%
}\hspace{-10pt}
\subfigure[risk-free reference point]{
        \begin{minipage}[t]{0.3\textwidth}
                \centering
                \includegraphics[width=2in]{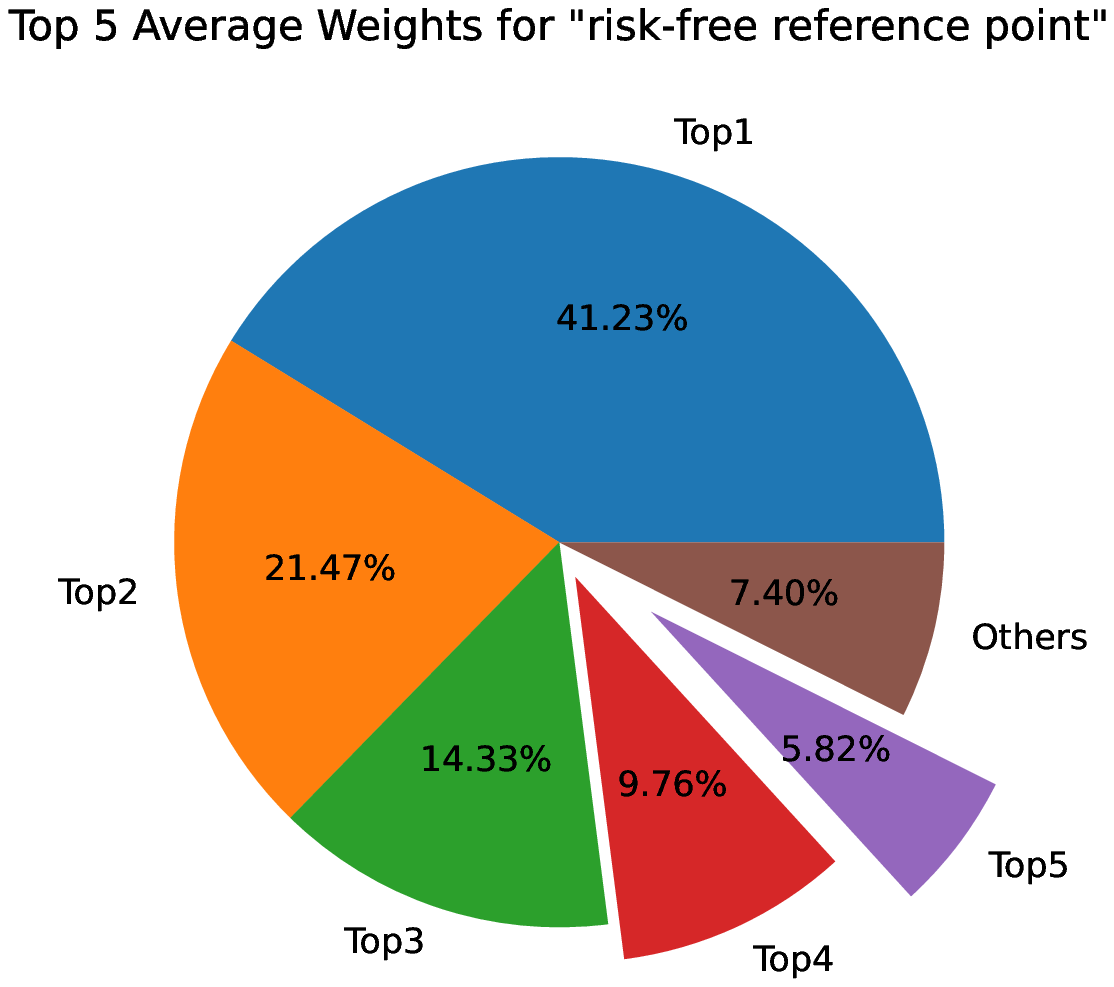}
        \end{minipage}%
}\hspace{-10pt}
\subfigure[large reference point]{
        \begin{minipage}[t]{0.3\textwidth}
                \centering
                \includegraphics[width=2in]{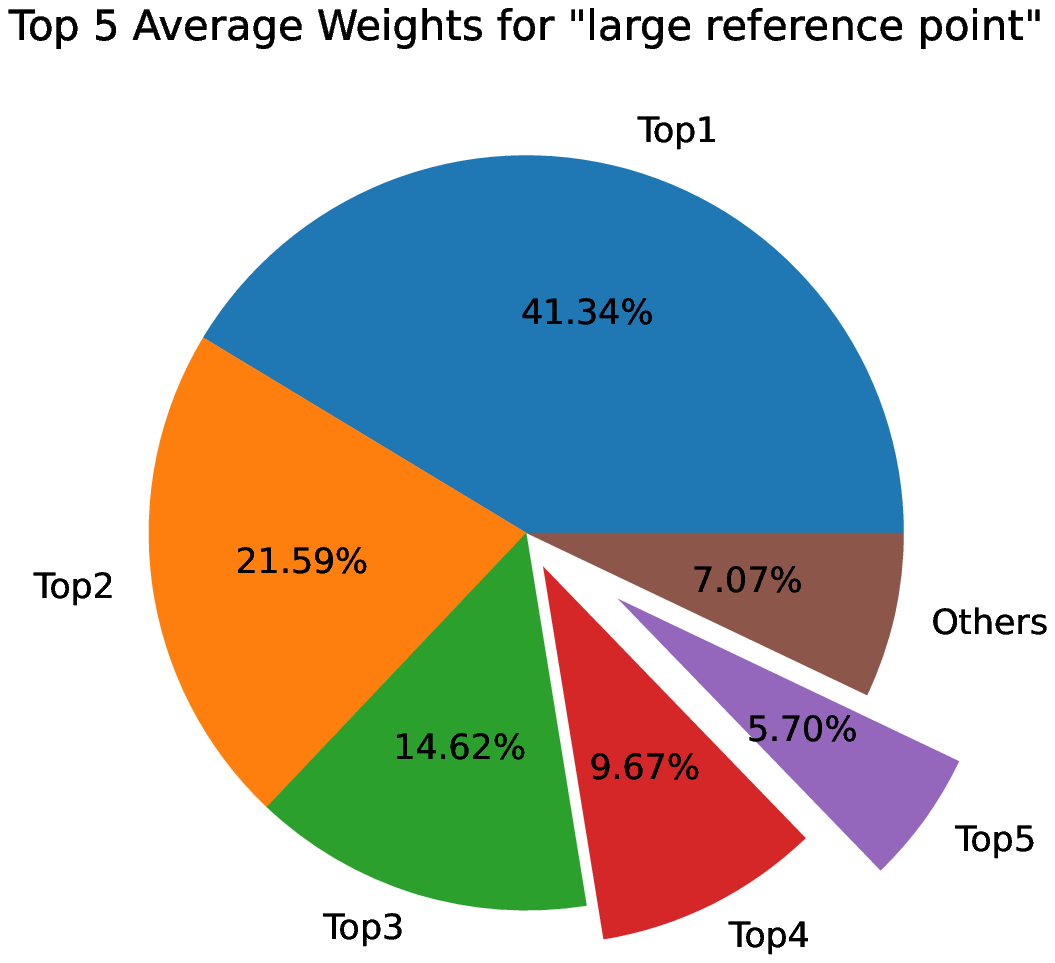}
        \end{minipage}%
}\hspace{-10pt}

\subfigure[no risk aversion and risk seeking]{
        \begin{minipage}[t]{0.3\textwidth}
                \centering
                \includegraphics[width=2in]{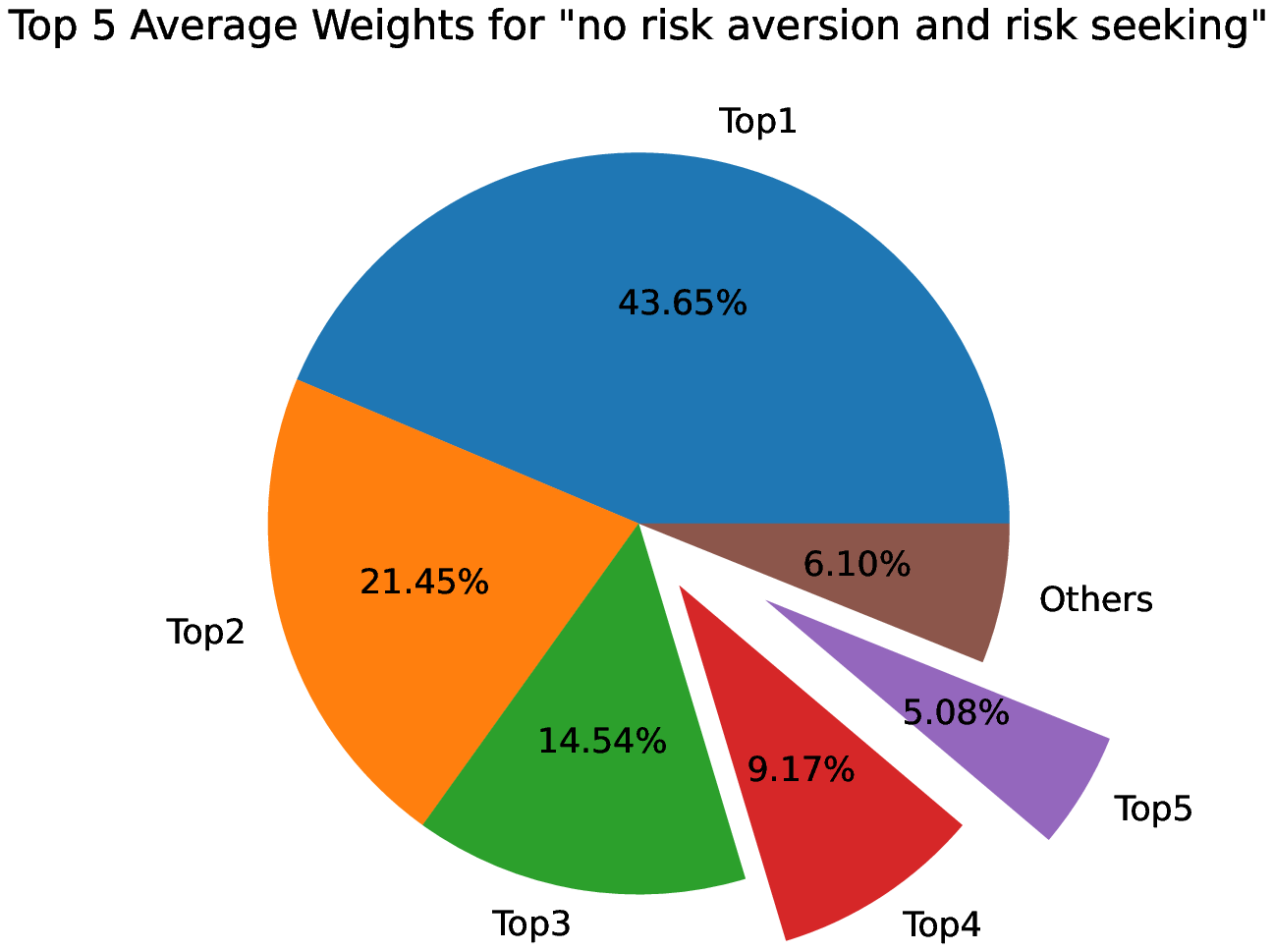}
        \end{minipage}
}\hspace{-10pt}
\subfigure[no loss aversion]{
        \begin{minipage}[t]{0.3\textwidth}
                \centering
                \includegraphics[width=2in]{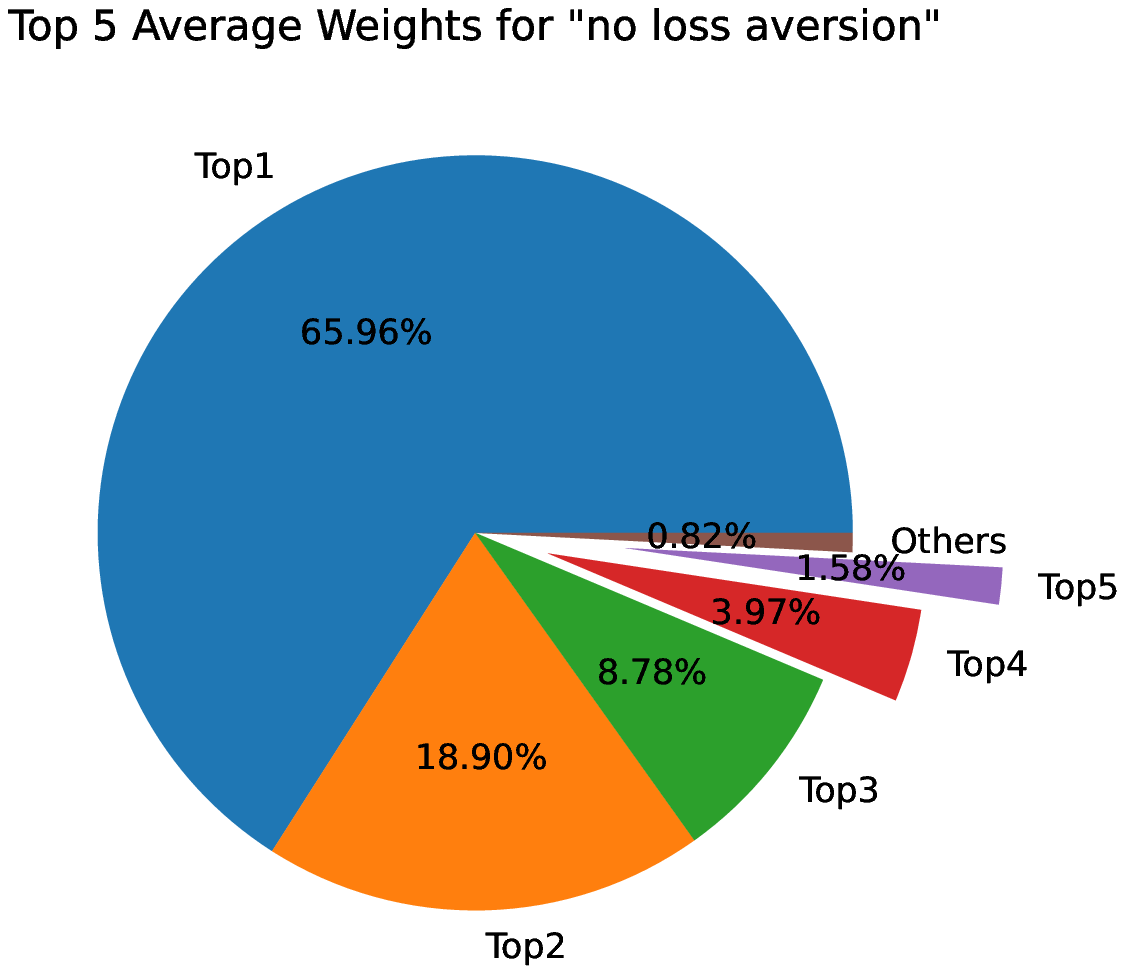}
        \end{minipage}
}\hspace{-10pt}
\subfigure[no probability distortion]{
        \begin{minipage}[t]{0.3\textwidth}
                \centering
                \includegraphics[width=2in]{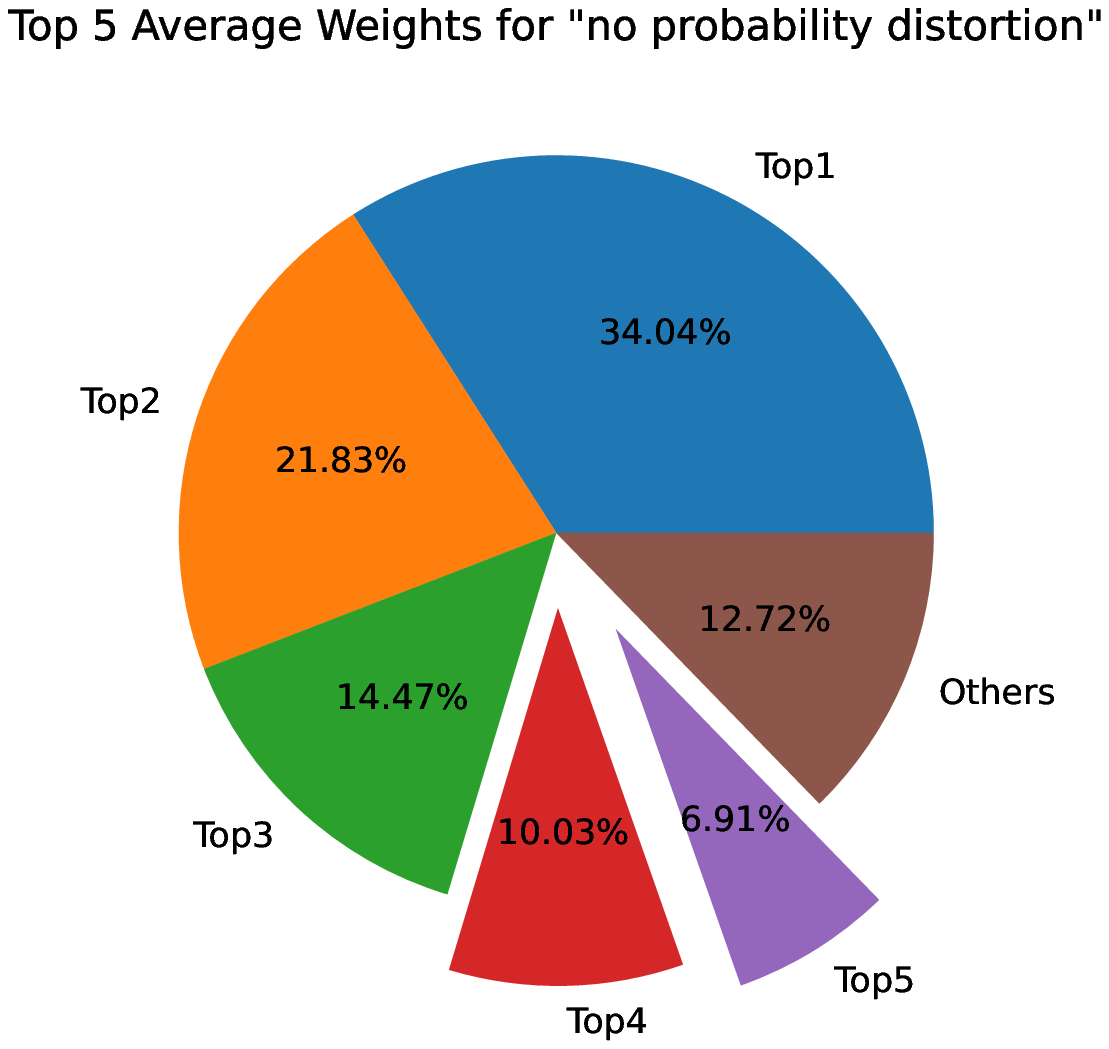}
        \end{minipage}
}\hspace{-10pt}
\centering
\caption{Average weights on five top and other assets.}
\label{fig:five}
\end{figure}

 Next, we explored the diversification levels of the six optimal portfolios by scrutinizing their portfolio weights. Each optimal portfolio was assessed with 1000 portfolio weight vectors during the out-of-sample period. From each weight vector, we extracted the top five values of the coordinates and aggregated the remaining 43 coordinates.
By averaging these values across the out-of-sample period, we derived the average weights of the portfolios allocated to the top five and other assets, as depicted in Figure \ref{fig:five}.
The weights in the benchmark and risk-free reference point models are quite similar. However, the large reference point model assigns slightly higher weights to the top one, two, and three assets, while the no risk aversion and risk-seeking model allocates larger weights to the top one and three assets. The absence of loss aversion results in an extreme allocation strategy, heavily concentrating investments in the top one asset. On the other hand, the no probability distortion model demonstrates superior diversification, indicated by its lower weight on the top one asset and higher weights on other assets. This finding aligns with the notion that probability distortion can impede diversification by encouraging investors to engage in risky gambling behaviors, as discussed by \cite{barberis2013}. In summary, our analysis highlights the significant impact of behavioral biases on portfolio diversification. It underscores the potential risks associated with deviating from rational decision-making in investment strategies.

\begin{figure}[H]
        \centering
        \subfigure[risk-free reference point]{
                \begin{minipage}[t]{0.45\linewidth}
                        \centering
                        \includegraphics[width=2in]{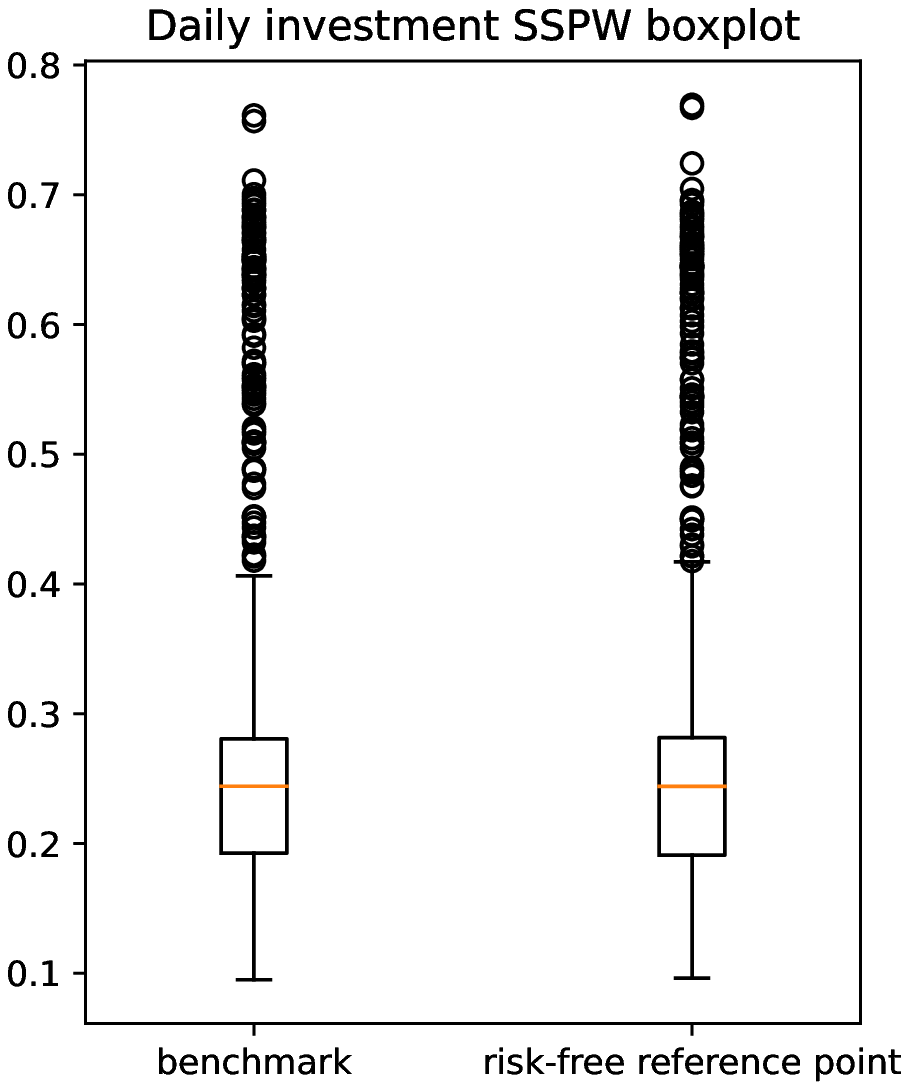}
                \end{minipage}%
        }\hspace{-10pt}
        \subfigure[large reference point]{
                \begin{minipage}[t]{0.45\linewidth}
                        \centering
                        \includegraphics[width=2in]{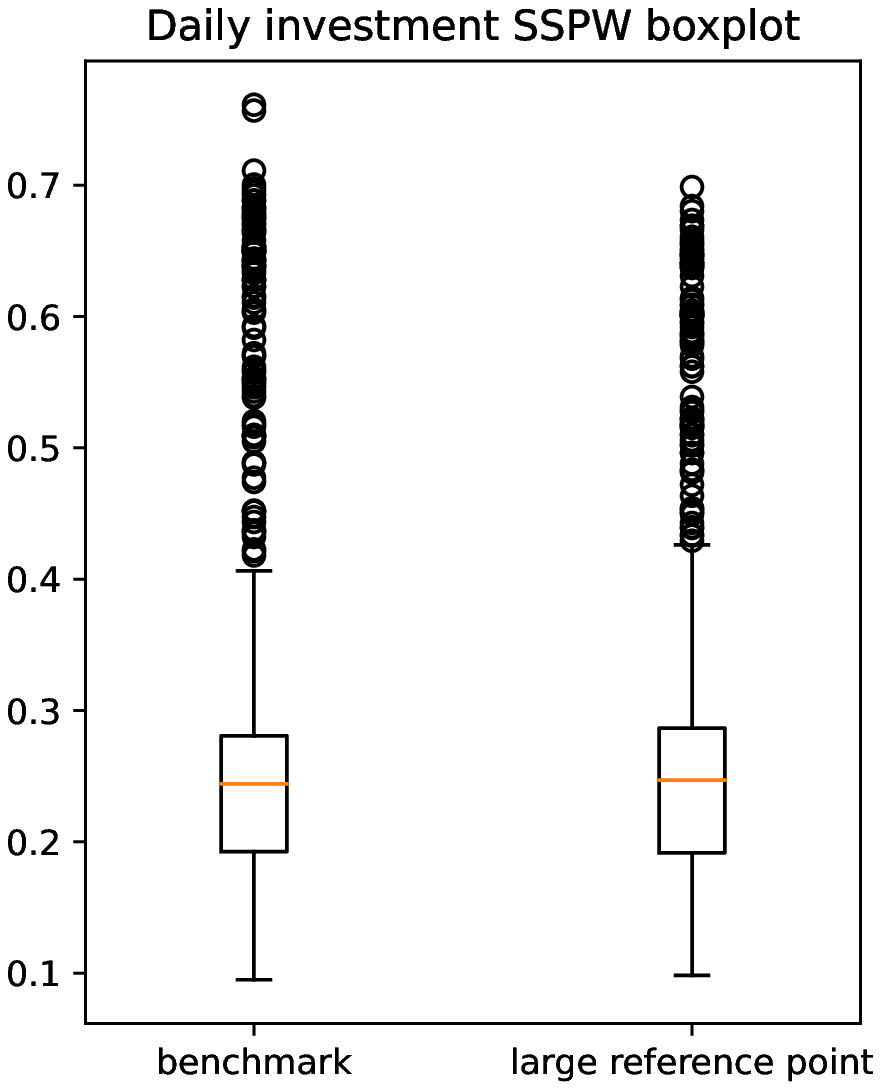}
                \end{minipage}%
        }\hspace{-10pt}

        \subfigure[no risk aversion and risk seeking]{
                \begin{minipage}[t]{0.3\linewidth}
                        \centering
                        \includegraphics[width=2in]{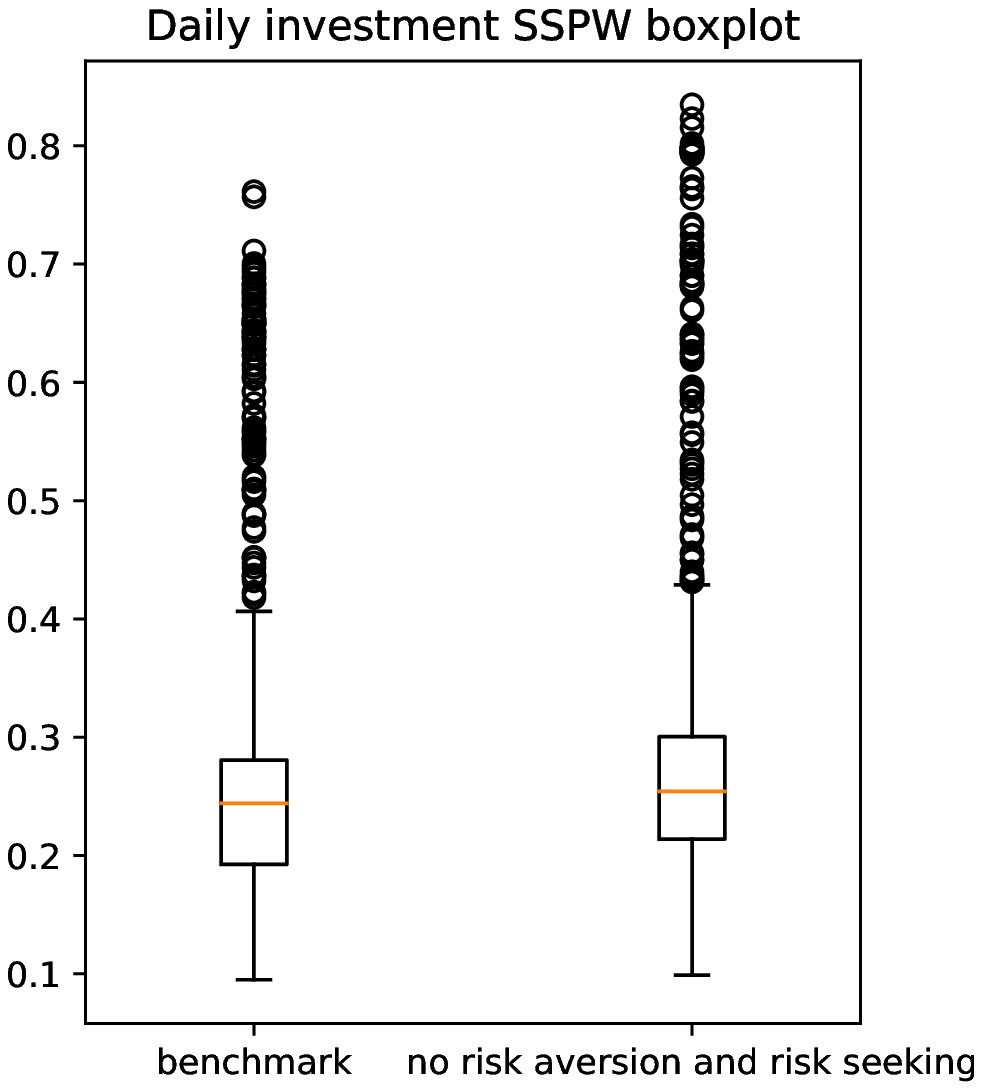}
                \end{minipage}%
        }\hspace{-10pt}
\subfigure[no loss aversion]{
        \begin{minipage}[t]{0.3\linewidth}
                \centering
                \includegraphics[width=2in]{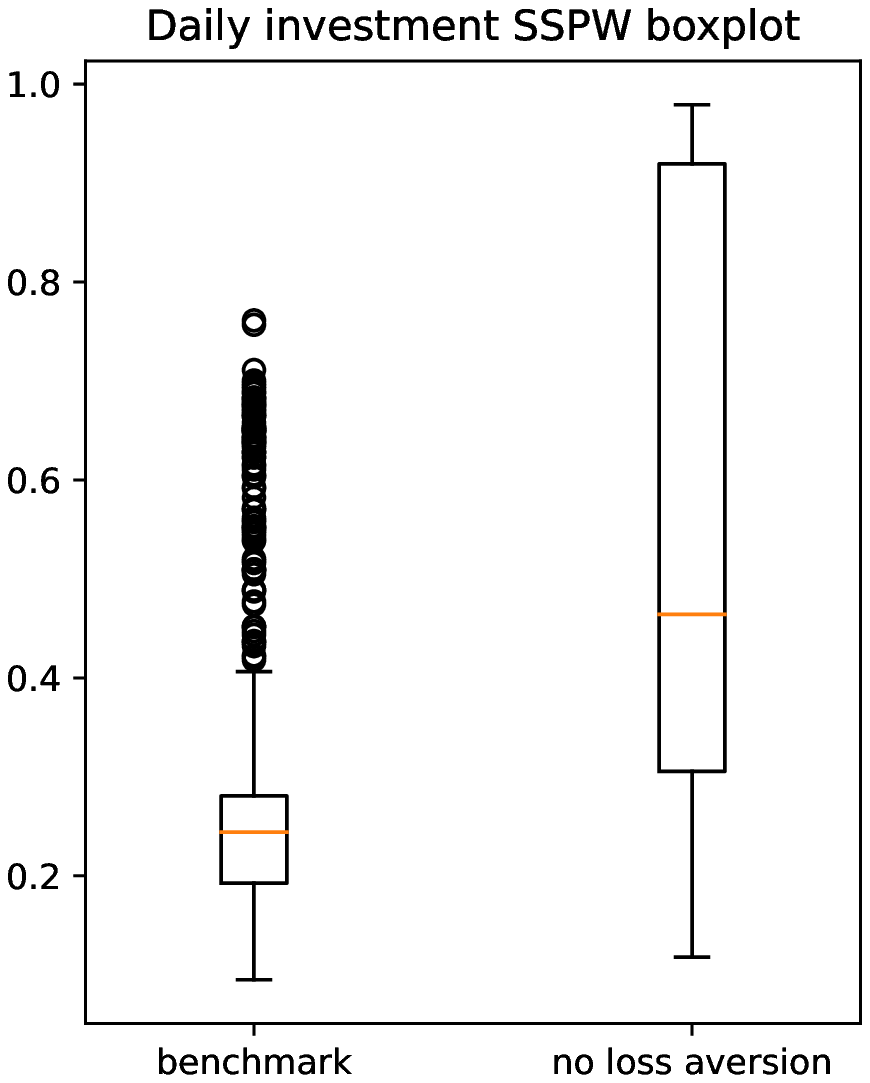}
        \end{minipage}%
 }\hspace{-10pt}
        \subfigure[no probability distortion]{
        \begin{minipage}[t]{0.3\linewidth}
                \centering
                \includegraphics[width=2in]{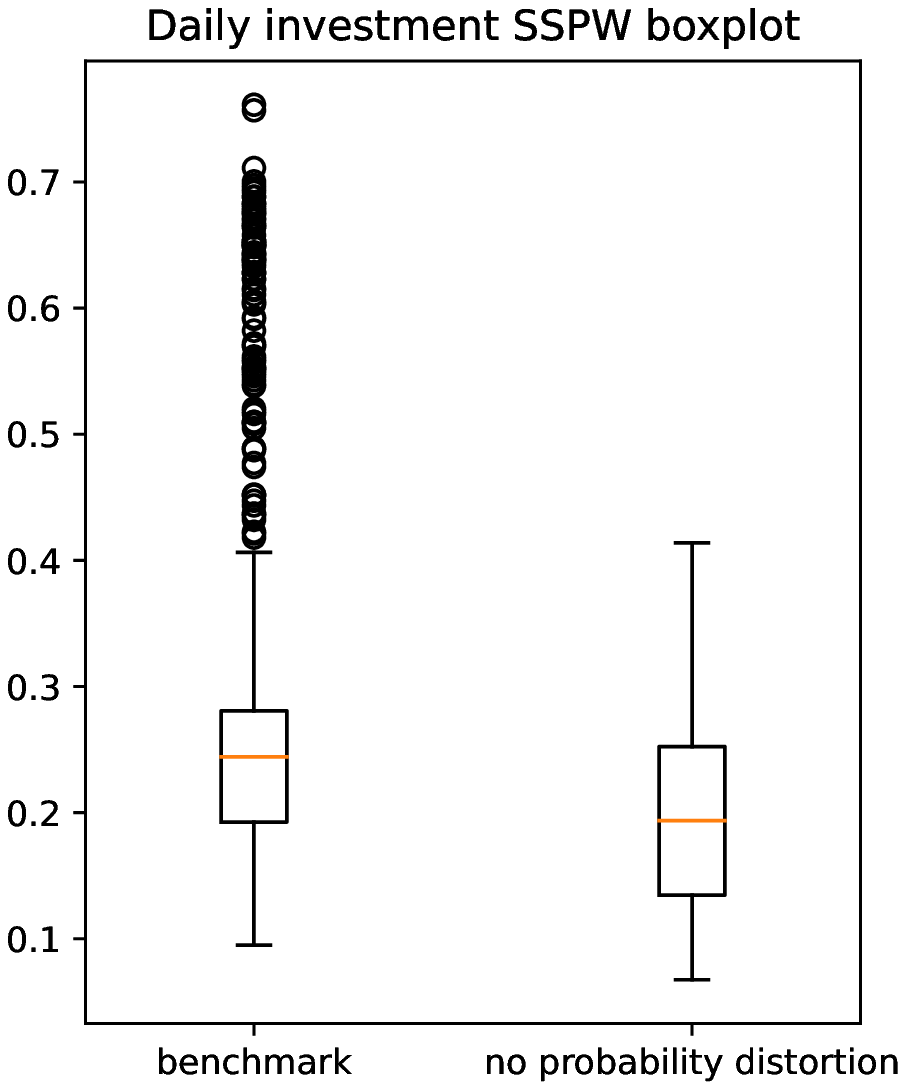}
        \end{minipage}%
}\hspace{-10pt}

        \centering
        \caption{The box-plots of six optimal portfolios' diversification degrees.}
        \label{fig:diversification}
\end{figure}

In addition to comparing the average portfolio weights of the different models, we also calculate the diversification degrees of the six optimal portfolios.
We follow \cite{goetzmann2008equity} and use the sum of squared portfolio weights (SSPW) as our metric, where SSPW is defined as follows:
\begin{align*}
SSPW = \sum_{i=1}^N\left(x_i - \frac{1}{N}\right)^2.
\end{align*}
Here, $N$ is the number of assets, and $x_i$ is the portfolio weight assigned to asset $i$ in the portfolio.
This metric measures how the portfolio weight $x_i$ deviates from the $1/N$ (equally weighted) portfolio, which is seen as the best-diversified portfolio. Hence, the lower the value of SSPW, the higher the level of diversification. We compute SSPW values of the six optimal portfolios every day in the out-of-sample period and show the corresponding box plots in Figure \ref{fig:diversification}.
Figure \ref{fig:diversification}(a) shows the comparison result between the benchmark case and the risk-free reference point case. The box plots exhibit striking similarity, with the exception of a few extreme outliers. We attribute this deviation to the risk-free rate's proximity to zero.  Figure \ref{fig:diversification}(b) shows the comparison result between the benchmark case and the large reference point case. Although the large reference point case has a similar median as the benchmark case, it has fewer extreme outliers. Figure \ref{fig:diversification}(c) shows the comparison result between the benchmark case and no risk aversion and risk seeking case. Removing risk aversion and risk seeking may result in a larger median and more extreme outliers. Figure \ref{fig:diversification}(d) shows the comparison result for the no loss aversion and benchmark case. As expected, removing loss aversion makes the investor gamble heavily and therefore results in portfolio concentration, i.e., a large increase in terms of SSPW.  Figure \ref{fig:diversification}(e) shows that compared with the benchmark case, removing probability distortion decreases the value of SSPW, thereby increasing the level of diversification. This finding is consistent with the theoretical prediction of \cite{barberis2013}.

\end{document}